\documentclass[11pt,a4paper]{article}
\usepackage[novbox]{pdfsync}
\usepackage{enumerate}
\usepackage{enumitem}
\usepackage{tabularx}
   \usepackage[margin=0.5in]{geometry}
\usepackage[utf8]{inputenc}
\usepackage{graphicx}
\graphicspath{{Figures/}}
\usepackage[T1]{fontenc}
\usepackage{fancyhdr, float}
\usepackage{mathrsfs,amssymb,amsmath,amsfonts,booktabs,array,xcolor,url,cite,
color,multirow,multicol}
\usepackage{slashbox}
\usepackage{stfloats}
\usepackage{extarrows}
\usepackage{textcomp}
\usepackage{ulem}
\usepackage{cases}
\usepackage{bm}
\usepackage{arydshln}
\usepackage{hyperref}\usepackage{array,subcaption}\usepackage{verbatim}
\usepackage[font=footnotesize]{caption}
\usepackage[all]{xy}
\usepackage{listings}
\usepackage{tikz, pgf,graphicx, pgfplots, color, xcolor}
	\pgfplotsset{compat=1.12} % indicates version of pgfplots used for consistency
\usepackage{tikz-3dplot}
\usepackage{tikz-cd}
	\usetikzlibrary{positioning}

\newtheorem{theorem}{Theorem}
\newtheorem{remark}{Remark}
\newtheorem{proposition}{Proposition}
\newtheorem{definition}{Definition}
\newtheorem{corollary}{Corollary}
\newtheorem{example}{Example}
\newtheorem{lemma}{Lemma}
\newtheorem{assumption}{Assumption}

 \def\eeD{\end{definition}} \def\beD{\begin{definition}}
\def\beR{\begin{remark}} \def\eeR{\end{remark}}
\def\beL{\begin{lemma}} \def\eeL{\end{lemma}}
\def\beC{\begin{corollary}
  }\def\eeC{\end{corollary}}
  \def\beT{\begin{theorem}}\def\eeT{\end{theorem}}
  \def\beP{\begin{proposition}} \def\eeP{\end{proposition}}
\def\beXa{\begin{example}} \def\eeXa{\end{example}}
\def\beA{\begin{assumption}} \def\eeA{\end{assumption}}
\def\im{\item}\def\com{compartment}     \def\vi{\overset{\rightarrow}{\i}}

   \def\ei{a_i}   \def\mV{{\mathcal V}}
  \def\Prf{{\bf  Proof:}}  \def\NGM{next generation matrix}
\def\brn{basic reproduction number}\def\DFE{disease-free equilibrium}
    \def\rd{\r_{dfe}} \def\sd{\s_{dfe}}

 \newcommand{\e}{\;\mathsf e}
    \def\va{\vec \al } 

\renewcommand{\i}{\;\mathsf i}
\renewcommand{\r}{\; \mathsf r}

\def\fe{for example} \def\wk{well-known}
\def\BEN{\begin{enumerate}}  \def\BI{\begin{itemize}}
\def\EEN{\end{enumerate}}   \def\EI{\end{itemize}}
\newcommand{\beq}{\begin{eqnarray}
    }
\def\eeq{\end{eqnarray}}
   \newcommand{\be}[1]{\begin{equation}\label{#1}}
\newcommand{\ee}{\end{equation}}
\def\bea{\begin{eqnarray*}} \def\im{\item} \def\Lra{\Longrightarrow}  \def\eqr{\eqref}  
\def\no{\nonumber} 

\newcommand{\R}{\mathbb{R}}

 \def\sd{\s_{dfe}}\def\al{\al}\def\sec{\section}
\def\b{\beta} \def\g{\gamma}    \def\de{\delta}
\def\z{\zeta}   
  
   \def\al{\al}
\def\L{\lambda} \def\La{\Lambda}\def\mR{\mathcal R} \def\fr{\frac}
\def\bc{\begin{cases}
  }
\def\ec{\end{cases}}
\def\bea{\begin{eqnarray*}}\def\ssec{\subsection} 
  
  \def\QED{\hfill {$\square$}\goodbreak \medskip} \def\satg{satisfying }  \def\sats{satisfies}  \def\saty{satisfy}
\def\eea{\end{eqnarray*}} \def\T{\widetilde}
   \def\I{\infty}
\def\no{\nonumber}\def\mF{\mathcal F}  
\long\def\symbolfootnote[#1]#2{
\begingroup
\def\thefootnote{\fnsymbol{footnote}}\footnote[#1]{#2}
\endgroup}
\def\fn{\symbolfootnote}
\providecommand{\pp}[1]{\left[#1\right]} %[.]
\providecommand{\pr}[1]{\left(#1\right)} %(.)
\def\bep{\begin{pmatrix}} \def\eep{\end{pmatrix}}
\def\bev{\begin{vmatrix}} \def\eev{\end{vmatrix}}

  \def\resp{respectively}

\newcommand\CRN{chemical reaction network }

\newcommand\CRNs{chemical reaction networks}

\newcommand{\bb}{\mathbf{b}}

\newcommand{\xb}{\mathbf{x}}
\newcommand{\yb}{\mathbf{y}}
\newcommand{\Yb}{\mathbf{Y}}
\newcommand{\Ab}{\mathbf{A}}

\newcommand{\Zb}{\mathbf{Z}}
\newcommand{\Sb}{\mathbf{S}}

\def\m0{{\mathcal R}_0} \def\la{\label}
\def\Eq{\Leftrightarrow}\def\qu{\quad}

\newcommand\s{\boldsymbol{s}}

\renewcommand{\epsilon}{\varepsilon}	% varepsilon
\newcommand\0{\mathbf{0}}

\renewcommand{\cdot}{}%The way we use it, we do not neet this dot

\newcommand*{\Scale}[2][4]{\scalebox{#1}{$#2$}}%
\def\frt{furthermore }

\newcommand{\inD}[1][\relax]{\def\argone{#1}\def\temprelax{\relax}
\ifx\argone\temprelax\right.\else\,\middle|#1\right.{}\fi}

\allowdisplaybreaks

\def\bb{\bff \beta} \newcommand{\bff}[1]{{\mbox{\boldmath$#1$}}}
\def\bzn{basic replacement number}  \def\ith{it holds that }
\def\s{\;\mathsf s}\def\al{\alpha}
\def\wl{w.l.o.g. } 
\renewcommand{\i}{\;\mathsf i}
\renewcommand{\r}{\; \mathsf r}
\def\m{b}\def\sec{\section}\def\bb{\bff \beta}\def\ME{mathematical epidemiology}
\def\AB{$(A,B)$ Arino-Brauer epidemic models}
\def\Dt{Descartes type}
\newcommand\red[1]{\textcolor{red}{#1}}
\newcommand\blue[1]{\textcolor{blue}{#1}}
\def\La{\Lambda}\def\m{\Lambda}\def\z{\;\mathsf z}
\newcommand{\figu}[3]{
\begin{figure}[H]
\centering
\includegraphics[scale=#3]{#1}
\caption{#2\label{f:#1}}
\end{figure}
}
\begin{document}
\title{ Advancing Mathematical  Epidemic Modeling via synergies with  Chemical Reaction Network Theory and Lagrange-Hamilton Geometry}
	\author{
 Florin Avram$^{1}$,  Rim Adenane$^{2}$, Mircea Neagu$^{3}$,
}
\maketitle

\begin{center}
	
		$^{1}$
		Laboratoire de Math\'{e}matiques Appliqu\'{e}es, Universit\'{e} de Pau, 64000, Pau,
 France; avramf3@gmail.com \\
$^{2}$ \quad Laboratoire des Equations aux dérivées partielles, Alg\'ebre et G\'{e}om\'{e}trie spectrales, d\'epartement des Math\'ematiques, Universit\'e Ibn-Tofail, 14000, Kenitra,
 Maroc; rim.adenane@uit.ac.ma \\
 $^{3}$ \quad  Transilvania University of Bra\c{s}ov, Department of Mathematics and Computer Science, 500091, Bra\c{s}ov, Romania; 	mircea.neagu@unitbv.ro
	\end{center}

\begin{abstract}

 This essay reviews some key concepts in   mathematical epidemiology and examines the intersection of this field with related scientific disciplines, such as  chemical reaction network theory and Lagrange-Hamilton geometry.
 Through a synthesis of theoretical insights and practical perspectives, we underscore the significance of essentially non-negative kinetic systems in the development and implementation of robust epidemiological models. Our purpose is to make the case that currently mathematical modeling of epidemiology  is focusing too much on simple particular  cases, and  maybe not enough on more complex models, whose  challenges would require cooperation  with scientific computing experts and with researchers in the ``sister disciplines" involving   essentially nonnegative kinetic systems (like virology, ecology,
\CRNs, population dynamics, etc).  %(as an example, as argued below,  even  the  SIR model is not yet sufficiently understood). A possible alternative to this focus on examples would be to search for general principles, which apply also to the other applied sister disciplines  .

\end{abstract}


\tableofcontents
{%\huge

\section{Introduction}
{\bf Motivation}. Dynamical systems  are  very important in the  ``sister sciences" of mathematical epidemiology (ME), virology, ecology, population dynamics, \CRN\  (CRN) theory, and other related domains. They may have a priori very complex behaviors, including in two dimensions (as illustrated by Hilbert 16'th problem, for example).
   %\fn[4]{For example, in  the  ecology literature on {\bf three dimensional food-chains},  chaos was noticed to appear first numerically,  and then explained via slow-fast expansions  in a long series of papers \cite{klebanoff1994chaos,kuznetsov1996remarks,kuznetsov2001belyakov,deng2001food,deng2003food,deng2004food,deng2006equilibriumizing,deng2017numerical}.}

Below, we will offer   speculate on questions like: what is a (``useful") mathematical epidemiology ODE model? is it different from a mathematical virology ODE model? from a \CRN\ ODE model? do these disciplines have general results, besides those which apply to all ODE dynamical systems? and if they do, are their results which apply
specifically to one of these ``sister sciences", and not to the others?

We will offer a specific, though highly disputable answer to the first question; as for the other questions, we will just speculate now that the answers are No, No, Yes, No.

In fact, the positive answer Yes is not a speculation. It has long been noted that the natural restriction to essentially nonnegative (i.e. positivity preserving)
 ``mass action kinetics" leads often to   results which are  surprisingly simple,
 despite high-dimensionality (see for example \cite{gunawardena2003chemical}).  This has motivated
 several researchers  in the applications cited above  to propose turning  mass-action polynomial kinetics into a unified tool for  studying all applied disciplines involving essentially nonnegative dynamical systems. One striking result about this class is its  equivalence  to the ``absence of negative cross-terms", has been discovered in \cite{hun}. Unfortunately, this result seems unknown outside \CRN\ theory, and it has been reproved in particular examples an uncountable number of times.

There are not many  sources which attempt to develop a unified view-point of essentially-nonnegative systems; for one exception, see   the book \cite{Haddad}. Further striking results are slow in coming. Some hope arises though from
 recent  papers which might be associated to  the new unifying banner of ``algebraic biology" -- see for example \cite{pachter2005algebraic,macauley2020case,torres2021symbolic}.

 Unfortunately, the opposite to unification is happening. The sister sciences seem to grow further and  further apart, due  to their focus on examples, to the point that a suspicious reader might ask himself whether sometimes they might not be studying the same system  under different names, without being aware of each others results.

 The challenge of mathematical epidemiology and the sister disciplines comes maybe from the fact that they must  split their efforts in at least 4 directions (besides the not to be forgotten identification of statistical data). These are the: \BEN   \im  creation of appropriate mathematical models,
 \im development of general results, \im creation of scientific computing tools for resolving classes of general models, \im  solving of  simple specific examples, usually via "hand computations". \EEN
 
  The last direction is by far the most represented in the literature. However, it is natural to hope that eventually the resolution of
examples will end up contributing to the ``unification" body of knowledge
 subsumed in the first three directions. 

 {\bf A crucial difference between \ME\ and \CRN\ theory}.  The most fundamental aspect  of mathematical epidemiology is the existence of at least two possible special fixed states.
   The first, the  DFE,  corresponds to the elimination of all compartments involving sickness. Note that the absence of boundary states is typically assumed in the nowadays \CRN\ literature, so that
 these two bodies of knowledge   diverge from the start, due to the examples they focus on.

 \beR \label{r:0} Mathematical epidemiology becomes easier if one restricts  to the subset of all  possible essentially non-negative dynamical systems which admit a unique DFE, like in \cite{Van08}. It may be that further restrictions could be beneficial for structuring the field.
For example, one could restrict as well to dynamical systems for which the unique DFE is locally stable for a non-void subset of parameter values, or to restrict even further to  particular classes, like that to be discussed in  section \ref{s:ker}. Unfortunately, fundamental  level $0$ questions like this seem very difficult to investigate.

Another  fundamental  level $0$ question, to be discussed more in detail below, is  the computation of  $R_0$. Defined nowadays as the spectral radius of the so called ``\NGM" (NGM)  \cite{Diek,Van,Van08}, this quantity must  
\BEN \im be such that  $R_0 <1$ defines the local stability domain of the unique DFE,   and 
\im reflect numerically the tendency of an  epidemics to spread itself, or serve as ``proxy for pathogen virulence" \cite{brouwer2019phenotypic}.
\EEN
Unfortunately, the NGM is obtained via a decomposition which was not specified uniquely in the  papers cited above. This leads to the non-uniqueness of $R_0$, and leaves its final choice to the latitude of the ``expert epidemiologist", a situation which is maybe not ideal.

We propose to complement the classical NGM method by a unique specification of   F  as the positive part of all the interactions which involve both  input/susceptible and infectious classes.

%even though an obvious ``natural" such specification, to be spelled out below, exists.
 We are not sure whether  there are  benefits in allowing for other NGM's besides the ``natural" one. 
Based on  experiments with such subclasses, we  offer as a potential definition of mathematical epidemiology models that of essentially non-negative models for which the {\bf natural NGM} associated to an unique DFE identifies correctly the local stability domain for the dynamical system, for a certain nonvoid domain of values of the parameters.
  \eeR

{\bf Contents}. Section \ref{s:ME} offers ``a bird's eye view" of the main result in mathematical epidemiology, the ``stability threshold theorem". This is accompanied in {Subsection} \ref{s:SEIR} by some Mathematica scripts which allow the reader to check our computations, in the particular case
of the   classic and easy to study four
\com s SEIR epidemic model.   %The point of this example  is to suggest a  reduction to a two-dimensional model without boundary. Unfortunately, our  reduction was too complicated to study for us (even though the original model is easy).

\iffalse
The bulk of our paper is dedicated to mathematical epidemiology and virology. In particular,  being especially aware of how far from achieving the four  goals listed in the introduction we are in  these fields,  we offer a critic in  Section \ref{s:crit} -- see also Remark \ref{r:0}.

Section \ref{s:Gro}  explains  how we succeeded to simplify    tedious hand-computations in papers like \cite{Jin}, and generalize them, via the use of Groebner bases, and the introduction of new concepts of ``Groebner eliminated determinants, traces, and Hurwitz determinants" (or, more generally, ``Groebner eliminated characteristic polynomial coefficients").  We believe these concepts deserve  further research.
\fi

Section \ref{s:back} presents some  general background information on  kinetic models, using the unifying notations of \CRN\ theory.

Section \ref{s:Neag} reviews some {Lagrange-Hamilton geometric objects which may be associated to any dynamical system, and the Subsection \ref{s:n1}
computes them for SIR-PH-FA models.} %We conjecture that these objects will turn out useful in slow-fast analysis of dynamical systems, and in particular in ME.

\iffalse

In section \ref{s:SIR}, we turn to our main example, an   eleven  parameter  three \com s SIR-type varying population epidemic model, which generalizes both a SIR model studied in \cite{AABH,NillI,NillII},  and a superinfection case studied in \cite{Mogh,Jin}, with the probable motivation of illustrating the possibility of Hopf bifurcations. This model may be of interest also to \CRN\ theorists, since it is not that  far from their ``auto-catalator" model. The presence or Hopf bifurcations
for our SIR model is still mostly uncharted theory, even in the simplest particular cases presented in Sections \ref{s:Nill},\ref{s:Jin}.

Section \ref{s:Hopf} offers a brief review of saddle-node, Poincaré-Andronov-Hopf and Bogdanov-Takens bifurcations. The subsection \ref{s:opt} includes the so called ``Bioswitch optimization" approach of \cite{otero2017chemical,yordanov2020bioswitch} for finding saddle-node bifurcations, whenever the analytic computations become too hard
  (we prefer the name ``Routh-Hurwitz optimization").   We also note that in principle it can  be easily modified to apply  to Hopf and Bogdanov-Takens bifurcations. This application of ideas originating from  biochemical researchers illustrates well the interest of a unification of the ``sister disciplines" which apply mathematical mathematical concepts to biology.

%Section \ref{s:bif} presents  bifurcation  diagrams for our examples.
\fi

Finally, we end with  a brief review of the field of chemical reaction networks (which is focused mostly  on the study of the stability and bifurcations of  fixed interior points)-- see   subsection \ref{s:CRN},   and represent in section \ref{s:3} several examples of epidemic models  as \CRNs, illustrating  the often mentioned, but  never exploited fact that
polynomial epidemic  models
are  usually   mass-action systems.
We are fascinated by this discipline, since CRN researchers succeeded in developing several general laws concerning their dynamical systems. It is natural to ask if these might turn useful in other disciplines, but  as far as we know,  this has not happened yet. This may possibly be due to the fact that boundary fixed points are ubiquitous  in epidemiology, immunology, population science, etc, but assumed not to exist
in \CRN\ theory.

\section{A bird's eye view of \ME: the \DFE, the \NGM, and an algorithmic definition  of a stability threshold associated to the \brn\ $R_0$ \label{s:ME}}

\ssec{The disease free equilibrium (DFE)}
  The DFE may be defined as a ``maximal boundary state", and    may be  found by identifying   a maximal sub-system  which factors \begin{equation}\label{Kfac}\vi'= M \vi.\end{equation}
  The components  $\vi$ will be called infectious states, and the set of its  indices will be denoted by ``infec".  Note that specifying ``infec"  induces a partition of both the coordinates and the equations into infectious (eliminable) components, and the others.

     The DFE is easily  computed by  solving the remaining
 ``non-infectious equations" with $\vi=0$.   In this paper we will assume its  uniqueness, at least after excluding   biologically irrelevant fixed points, like an unreachable origin. 
 
 \iffalse 
 
  We give now a      very elementary script, to emphasize the fact that any ODE model ``mod" (like SIR, etc...), is   a pair mod= (dyn,X) consisting of a vector field and a list of variables. But, since sometimes only numeric solutions are possible, our fixed point Mathematica script below has also an optional numerical condition parameter ``cn".

  \begin{verbatim}
 DFE[mod_,inf_,cn_:{}]:=Module[{dyn,X},
  dyn=mod[[1]]/.cn;X=mod[[2]];
  Solve[Thread[dyn==0]/.Thread[X[[inf]]->0],X]];
 \end{verbatim}

 \begin{verbatim}
 FP[mod_,cn_:{}]:=Module[{dyn,X},
   dyn=mod[[1]];X=mod[[2]];Solve[Thread[(dyn//.cn)==0],X]];
 \end{verbatim}
 \fi

 %\input{R0n}

\beR

There are (at least) two flavors of \ME\ and two corresponding  formulas for $\mR_0$:\BEN \im One, for ODE/Markovian models, identifies $\mR_0$  as the spectral radius of the Perron-Frobenius eigenvalue of the ``NGM"  $F V^{-1}$, obtained by splitting  the infectious equations as $$\vi'=\vi (F-V),$$ where $F$ has only nonnegative elements, and $-V$ is a Markovian generating matrix (this result requires that the set ``infec" of infectious equations satisfies appropriate conditions -- see \cite{Diek,Van,Van08}, and guarantees neither existence of ``infec", nor its uniqueness).

\im Under the second, ``non-Markovian/renewal"  approach,  $\mR_0$ is computed as the integral of  an ``age of infection kernel"  \cite{Diek}.
\EEN

The intersection of these two classes, the ODE/Markovian and the non-Markovian/renewal models, is  the notable context of ``rank one SIR-PH-FA epidemic models" \cite{AABGH}, which are a particular case of the more general \AB\ introduced in \cite{AABBGH}. Alternatively, these are  precisely the renewal models with a matrix-exponential kernel.
 %this matrix is $B V^{-1}$, where $B$ is the matrix with nonnegative elements appearing in the ``new infections" part of the infectious equations, and $-V$ is a constant Markovian subgenerator matrix appearing  in  the infectious equations.
The equivalence of the two approaches for this class of epidemic models is proved very concisely in \cite{AABGH}, and it may also be read between the lines of the wider scope papers \cite{Diek18,Diek22}. We revisit this class in the next section

 Beyond this class of simple models, \ME\ is still largely a collection
 of examples and open problems. Even the classic 3 compartments SIR process has been fully analyzed only recently in \cite{NillI,NillII}. % We review this example below, and in particular the fact that the SIR may have two endemic points when $\mR_0>1$, when the state space decomposes into disconnected invariant sets, which corrects a conjecture made in \cite{AABH,AAH}.
\eeR
\subsection{The subclass of ODE epidemic models with one susceptible class and \NGM\  of rank one;  their Markovian semi-groups, age of infection kernels, and   $\mR_0$  formula   \la{s:ker}}
The idea behind the \NGM\ method is that the infectious components may be expressed in function of the others. This is especially easy to state for  {\bf SIR-PH-FA models} \cite{AABBGH},  defined by:
\begin{align}
\label{SYRPH}
\vi'(t) &=  \vi (t) \pp{\s(t)  \; B+ A
- Diag\pr{\bff{\de}+\La \bff 1}}:=
 \vi (t) \boxed{(-  V +\s(t) B)}  \no\\
 \s'(t) &= \pp{\La -\pr{ \La + \g_s } \s(t)}- \s(t) \T {i}(t)+ \g_r {\mathsf r}(t), \; \; \T {i}(t) = \vi(t) \bb \nonumber\\
 \bb&=\bep\bb_1\\\vdots\\\bb_n\eep, \quad \bb_i=(B \bff 1)_i=\sum_{j} B_{i,j}, i=1,...,n\nonumber\\
{\mathsf r}'(t) &=  \vi (t) \bff a+ \s(t)  \g_s   - (\g_r+\La){\mathsf r}(t), \qu \bff a=(-A) \bff 1.
 \end{align}
 Here,
 %the transmission rate $\b$ is a constant and
\BEN
\im $\s(t) \in \mathbb{R_+}$ represents the set of individuals susceptible to be infected (the beginning state).
\im ${\mathsf r}(t) \in \mathbb{R_+} $  models
recovered  individuals (the end state).
\im $\g_r$ gives  {the rate at which recovered individuals lose immunity,}
and  $\g_s$ gives the rate at which individuals are vaccinated (immunized). These two transfers connect directly the beginning and end states (or classes).
\im the row vector $\vi(t) \in \mathbb{R}^n$  represents the set of individuals in different disease states.
\im $\Lambda >0$ is the per individual death rate, and it equals also the global birth rate (this is due to the fact that this is a model for proportions).

\im  $A$ is a $n\times n$ {\bf Markovian sub-generator matrix} which describes transfers between the disease classes. Recall that a Markovian sub-generator
matrix for which each
off-diagonal entry $A_{i,j}$, $i\neq j$, satisfies $A_{i,j}\geq 0$, and such that the row-sums are non-positive,  with  at least one inequality being strict.\fn[4]{Alternatively, $-A$ is a  non-singular M-matrix \cite{Arino}, i.e.
  a real matrix $V$ with  $ v_{ij} \leq 0, \forall i \neq j,$ and having eigenvalues
whose real parts are nonnegative \cite{plemmons1977m}.}

The fact a Markovian sub-generator appears  in  our ``disease equations" suggests that certain probabilistic concepts intervene in our deterministic models, and this is indeed the case--see below. %Note also that typical epidemic models \saty\ $A_{i,j} A_{j,i}=0$, $i\neq j$, and so this matrix may be arranged  to be triangular.

\im $\bff \de \in \mathbb{R_+}^n$ is a column vectors giving the death rates caused by the epidemic in the  disease \com s. The matrix $-V$, which combines $A$  and the birth and death rates $\La, \bff \de$, is also a Markovian sub-generator.

\im $ B $ is a $n \times n$ matrix. {We will denote by $\bb$ the vector containing the sum of the entries in each row
of $B$, namely, $\bb= B \bff1$.} Its components $ \bb_i$ represent the
{\bf total force of infection} of the  {disease} class $i$, and $\s(t)   \vi(t) \bb$ represents the
total flux which must leave class $\s$.
Finally, each  {entry} $B_{i,j}$, multiplied by $\s$,  represents the
force of infection from the  {disease} class $i$ onto class $j$, and our essential assumption below will be that $B_{i,j}=\beta_i \al_j,$ i.e. that all forces of infection are distributed among the infected classes conforming to the same probability vector $\va=(\al_1,\al_2,...,\al_n)$.

\EEN

\beR  The Jacobian of the SIR-PH-FA model with $\g_r=0$ is:
\begin{equation}\label{JacPh} {\bep \s B-V& \vi B \\ -\s \bb& -\La -\g_s- \vi \bb \eep}\end{equation}
where $\bb$ is defined in \eqr{SYRPH}.
\eeR
\beR \label{r:PF}    Note that the factorization of the  equation for the diseased \com s $\vi$  implies a representation of $\vi$ in terms of $\s$:
\begin{equation}\label{irep}
\vi(t)=\vi(0) e^{- t V + B \int_0^t \s(\tau) d \tau} =\vi(0) e^{\pp{- t I_n + B V^{-1} \int_0^t \s(\tau) d \tau} V}. \end{equation}

{In this representation intervenes an essential character of our story, the matrix $B V^{-1},$ which is proportional to the \NGM\ $\s B V^{-1}$}. { A  second representation \eqr{vc} below will allow us to embed our models in the interesting class of distributed delay/renewal models, in the case when
$B$ has rank one.}

\eeR

\beP\label{p:ren}  Consider a  SIR-PH-FA model \eqr{SYRPH} with one susceptible class, without  $\g_r=0$, so that $\mathsf{r}(t)$ does not affect the rest of the system, and with $B= \bb \va$ of rank one. {Recall the total force of infection }
$$\T {i}(t) =  \vi(t) \bb.$$
Then
\BEN
\im  The solutions of the ODE system
\eqr{SYRPH}    \saty\ also a ``distributed
delay SI system"  of two scalar equations
\begin{align}\label{SI}
\bc \s'(t)=\Lambda -\pr{ \La + \g_s } \s(t)-\s(t) \T i(t)\\
\T i(t)=\vi(0)  e^{-tV} \bb+
\int_0^t s(\tau) \T i(\tau)
a(t-\tau) d\tau,\ec
\end{align}

where \begin{equation}\label{renker} a(\tau)=  \va e^{-\tau V} \bb, %\; \; -V=\bep  {-(\g_e +\La)} & \g_e %\\0 & -\pr{\g   +\La +\nu}\eep
 %\\\va e^{\tau A} \bb&\La=0\ec
 \end{equation}
with $-V=A
- \pr{Diag\pp{\bff{\de}+\Lambda \bff 1}}$  (it may be checked that this fits the formula on page 3 of \cite{Breda} for SEIR when $\Lambda= 0, \de=0$).\fn[3]{$a(t)$ is called ``age of infection/renewal  kernel; see  \cite{Hees,Brauer05,Breda,DiekHeesBrit,Champredon,Diek18,Diek22} for expositions of this concept.}
% since $$e^{\tau A}=\bep
 %e^{-\tau \gamma _e} & \frac{\gamma _e}{\gamma _e-\gamma }
 %(e^{-\tau \gamma  }-e^{-\tau \gamma _e} )\\
% 0 & e^{-\tau \gamma }\eep.$$

\im The \bzn\ $\mR$ has an integral representation
\begin{equation}\label{R0I}
\mR=\int_0^\I a(\tau)  d \tau= \int_0^\I \va e^{- \tau V} \bb d\tau =\va\ V^{-1} \; \bb.
\end{equation}

\EEN
\eeP

\Prf 1. The non-homogeneous   infectious equations may be transformed into an integral equation by applying the variation of constants formula.
The first step is the solution of the homogeneous part.  Denoting this by  $\Gamma(t)$,  \ith \begin{align}
\label{SYRIGI}
\vec \Gamma'(t) &=
- \vec \Gamma(t) V \Longrightarrow \vec \Gamma(t) = \vec \Gamma(0)e^{ t (-V)}.\fn[4]{When $\vec \Gamma(0)$ is a probability vector, \eqr{SYRIGI} has the interesting probabilistic interpretation of the survival probabilities in the various components  of the semigroup generated by the Metzler/Markovian sub-generator matrix $-V$ (which inherits this property from the phase-type generator $A$). Practically, $\vec \Gamma(t)$ will give  the expected fractions of individuals who are still in each compartment at time $t$.}
\end{align}

The variation of constants formula implies then that $\vi(t)$ satisfies \ the integral equation:
\begin{equation}\label{vc} \vi(t) = \boxed{\vi(0) e^{-t V} +
\int_0^t \s(\tau) \vi(\tau)  B e^{-(t-\tau) V} d \tau}.\end{equation}

Now in the rank one case $B=\bb \va$, and
 \eqr{vc} becomes
\begin{equation}\label{vcb} \vi(t) = \vi(0) e^{-t V}+
\int_0^t \s(\tau) \pp{\vi(\tau)  \bb} \va e^{-(t-\tau) V} d \tau.%:=\vi(0) e^{-t V}+
%\int_0^t \s(\tau) \vi(\tau)  \bb a(t-\tau) d \tau.
\end{equation}

Finally, multiplying both sides on the right by $\bb$ yields the result.

2. By the ``survival method"\fn[3]{This is a first-principles method, whose rich history  is
described in   \cite{Hees,Diek10}-- see  also \cite[(2.3)]{Champredon}, \cite[(5.9)]{Diek18}.}, $\mR$
 may  be obtained by integrating $\Gamma(t)$ with $\Gamma(0)=\va$.
 A direct proof is also possible by noting that all eigenvalues of the \NGM\ except one are $0$ \cite{Arino,AABBGH}.
\QED

\beR
The fact that DD systems can be approximated by ODE systems, {by approximating the delay distribution via one of Erlang, and more generally, of matrix-exponential type}, has long been exploited in the epidemic litterature, under the name of "linear chain trick" (which has roots in the Erlangization of queueing theory)-- see \fe\ \cite{Wearing2005,Feng,Wang2017,Diek18,cassidy2018recipe,Hurtado19,Ando,Diek22} for  recent contributions and further references.    The opposite direction however, i.e. identifying the kernels associated to
ODE epidemic models, seems  to have been forgotten.

\eeR

\subsection{Two examples: The SAIR/S$I^2$R/SEIR-FA epidemic model and the SLAIR/SEAIR epidemic model\label{s:SEIR}}

\iffalse Our ambition is to create a library of programs
  (in Mathematica, for example), which may be used  for ``solving any easy model", independently of the application it is used for.
   The first  script, ``JTD",  dedicated to outputing  Jacobians, traces and determinants
  \begin{verbatim}
 JTD[mod_]:=Module[{dyn,X,jac,tr,det},
    dyn=mod[[1]];X=mod[[2]];
    jac=Grad[dyn,X]/.cn;
    tr=Tr[jac];det=Det[jac];{jac,tr,det}];
 \end{verbatim}
 is both elementary and universally useful. Beyond this, we may need to take into account the special structure of mathematical epidemiology %(and virology, and population dynamics),
 which centers around the {\bf existence of boundary fixed points}.
 \fi
 \beXa \la{ex:SEIR}
The {9} parameters SAIR/S$I^2$R/SEIR-FA epidemic model \cite{Van08,RobSti,Ansumali,Ott,AAH},  its \NGM, and its \brn\ $\mR_0$:
\eeXa

 \figu{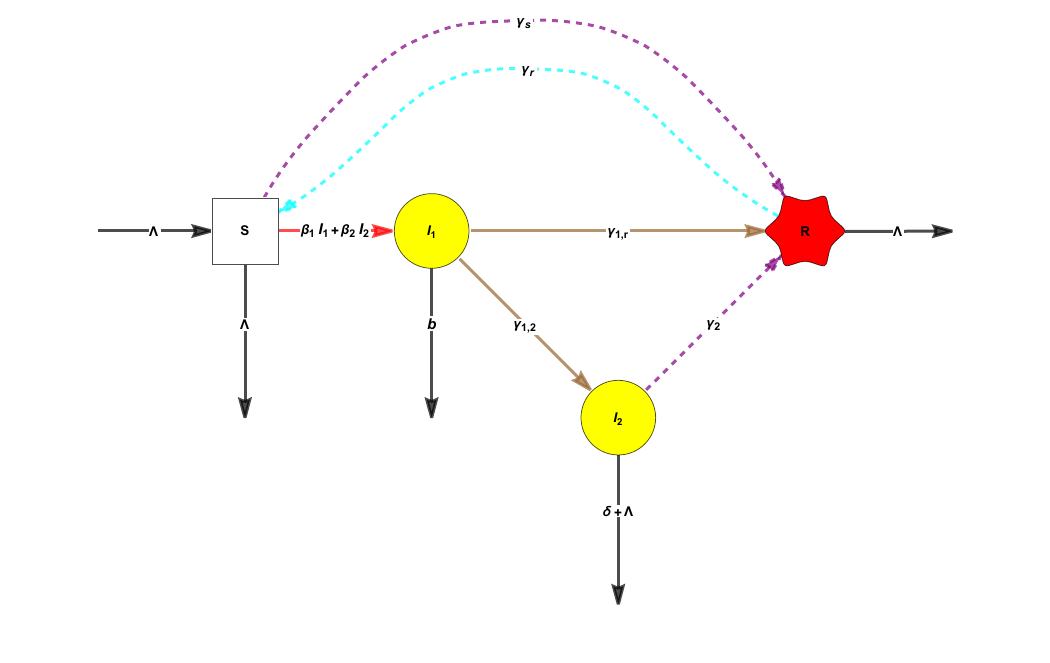}{Chart flow of the SI$^2$R model \eqr{SEIRsc}.  The red edge corresponds
to the entrance of susceptibles into the disease classes, the brown edges are the rate of the transition matrix V, and the cyan
dashed lines correspond to the rate of loss of immunity.  The remaining black lines correspond to the
inputs and outputs of the birth and natural death rates, respectively, which are equal in this case.}{.8}

 The proportions of  the {\bf long term, varying population model with $\de>0$}
 defined  in Figure \ref{f:SAIRS}   \saty\ {\bf approximatively} \cite{AABH,AAH}:

\begin{equation}\label{SEIRsc} \bc
\s'(t)= \Lambda   -\s(t)\pr{\beta_i  \mathsf i(t)+\beta_a \mathsf a(t)+\g_s+\Lambda} + \g_r \mathsf r(t)\\
\bep \mathsf a'(t)\\\mathsf i'(t)\eep =
\pp{\s(t) \bep  \beta_a     &  {\beta_i} \\
 0&  0 \eep- \bep    \g_e +\Lambda  &  0 \\
 -\ei& {\Lambda + \g_i + \de}\eep}
  \bep \mathsf a(t)\\\mathsf i(t)\eep
\\
\mathsf r'(t)=  \g_s \s(t)+ a_r   \mathsf a(t)+ {\g_i}   \mathsf i(t)- (\g_r+\Lambda) \mathsf r(t) \ec.
\end{equation}

 \beR \BEN \im  The  SAIR model is obtained when $a_r(\g_{1,r})=0=\de$ and the classic SEIR model is obtained when \frt\ $\beta_{a}=0$.
 \im We have written the ``infectious" middle equations to emphasize their factorization.
Also, for the factor appearing in these equations, we have emphasized a form
 \begin{equation}\label{V}F - V.\end{equation}
 Such decompositions, supposedly non-unique, not always existing \cite{Diek,Van,Van08}, are  used in the computation of the  next generation matrix (NGM)
 $$K=F.V^{-1}$$
 and of the \brn\ $R_0$.
\EEN \eeR

% \la{ex:SEIR}
 The SAIR/SEIR model
 is an \AB\ with  parameters $\va =\bep 1&0\eep, A=\bep -\g_1  &\g_{1,2}\\ 0&-\g_2   \eep, \bff a= (-A)\bff 1
 =\bep \g_{1,r} \\\g_2   \eep$ and
 \begin{equation}\label{sairsp}
  \bb=\bep \beta_ 1 \\ \beta_ 2\eep, \; \mbox{so}\; B=\bep \beta_ 1 & 0\\ \beta_ 2 & 0 \eep , \bff \de=\bep 0\\\de\eep, V=\bep    \g_1   +\Lambda  &  -\g_{1,2} \\
 0&  \g_2 +\Lambda +\de \eep.\end{equation}
The Laplace transform of the age of infection kernel  is:
 \begin{align} \label{aTr} \Hat{a}(s)&= \va (\s I+V)^{-1} \bb
 =
 \beta_ 1\frac{1 }{ \left(b+\g_1  +s\right)}+
 \beta_ 2 \frac{\g_{1,2}}{(b+\g_2  +\delta+\s ) \left(b+\g_1  +\s\right)},
 \end{align}
 %\frac{\g_{1,2} \beta_ 2 e^{-t (b+\g_2  +\delta )}-e^{-t \left(b+\g_1  \right)} \left(\beta_ 1 \left(\g_2  +\delta -\g_1  \right)+\g_{1,2} \beta_ 2\right)}{-\g_2  -\delta +\g_1  },$
 and the  Arino \& al. formula yields
$
\mR= \int_0^{\I} a(\tau) d \tau= \frac{\beta_ 1 (b+\g_2  +\delta )+\g_{1,2} \beta_ 2}{(b+\g_2  +\delta ) \left(b+\g_1  \right)}.$

The SIR model
 is also an  \AB\ with  parameters $\va =\bep 1\eep, A=\bep -\g   \eep,$
  $\bb= B= \beta, V=\bep    \g   +\La \eep.$

%\input{SEIRmat}
%\subsection{A generalized SLAIR epidemic model}
\beXa The SLAIR epidemic model \cite{YangBrauer,arino2020simple,AAGH} is defined by:
\begin{equation}\label{SLIARG}
\Scale[0.9]{ \bc
\s'(t)=   \Lambda -\s(t)\pr{\beta_ 2  i_2(t)+\beta_3 i_3(t)+\Lambda}\\
\bep \i_1'(t)& i_2'(t)& i_3'(t)\eep = \bep i_1(t)&i_2(t)& i_3(t)\eep
\pp{\s(t) \bep  0    &  0& 0 \\
 {\beta_ 2}&  0 &0\\
 \beta_ 3&0&0\eep+ \left(
\begin{array}{ccc}
 -\g_1 -\Lambda & \g_{1,2} & \g_{1,3} \\
 0 & -\g_2 -\Lambda  & \g_{2,3} \\
 0 & 0 & -\g_3 -\Lambda \\
\end{array}
\right)}
\\
{\mathsf r}'(t)= \gamma_{2,r}  i_2(t)+\gamma_3  i_3(t) -\Lambda \mathsf{r}(t)\ec.}
\end{equation}
This is an  \AB\ with  parameters
 \bea
 \Scale[0.9]{ \va =\bep 1&0&0\eep, A= \left(
\begin{array}{ccc}
 -\g_1  & \g_{1,2} & \g_{1,3} \\
 0 & -\g_2   & \g_{2,3}\\
 0 & 0 & -\g_3   \\
\end{array}
\right), \bff a= (-A) \bff 1
 =\bep 0\\ \g_{2,r}  \\\g_3  \eep, \bb=\bep 0\\ \beta_2 \\ \beta_ 3\eep, \;   \mbox{so}\; B=\left(
\begin{array}{ccc}
 0 & 0 & 0 \\
 \beta _2 & 0 & 0 \\
 \beta _3 & 0 & 0 \\
\end{array}
\right).}
\eea
The Laplace transform of the age of infection kernel  is:
 \bea
\Scale[1.1]{ \Hat{a}(s)=\beta _2 \frac{ \gamma _{\text{1,2}}}{\left(b+\gamma _1+\s\right) \left(b+\gamma _2+\s\right)}+ \beta_ 3 \pr{\frac{\g_{1,3}}{\left(b+\gamma _1+\s\right) \left(b+\gamma _3+\s\right)}+\frac{\gamma _{\text{1,2}} \gamma _{\text{2,3}}}{\left(b+\gamma _1+\s\right) \left(b+\gamma _2+\s\right) \left(b+\gamma _3+\s\right)}},}\eea
and the Arino \& al. formula yields
$
\mR= \frac{\beta _3 \gamma _{\text{1,2}} \gamma _{\text{2,3}}+b \beta _2 \gamma _{\text{1,2}}+\beta _2 \gamma _3 \gamma _{\text{1,2}}+b \beta _3 \gamma _{\text{1,3}}+\beta _3 \gamma _2 \gamma _{\text{1,3}}}{\left(b+\gamma _1\right) \left(b+\gamma _2\right) \left(b+\gamma _3\right)}.$
\figu{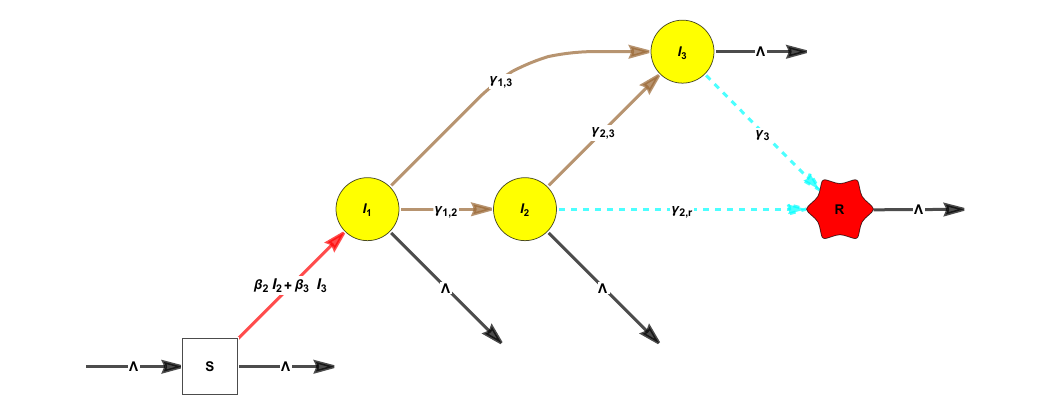}{ Chart flow of the SLAIR model \eqr{SLIARG}.
\iffalse The red edge corresponds
to the entrance of susceptibles into the disease classes, the dashed purple edges correspond to contacts
between infectious to susceptibles, the brown edges are the rate of the transition matrix V, and the cyan
dashed lines correspond to the rate of recovery.\fi}{.9} 
\eeXa

\subsection{A simple  recipe for $F-V$ gradient decomposition, and  its next generation matrix}

We complement here the famous $\mF-\mV$ ``equations decomposition" and \NGM\ method of \cite{Diek,Van,Van08}  by an algorithmic $ F-V$ decomposition.
\BEN \im
The user supplies the model ``mod" (a pair containing the RHS of the dynamical system, and its variables), and  the indices ``inf" of the disease (or infectious) variables; the indices of the
other compartments are denoted by ``infc".
\im Subsequently, the Jacobian of the infectious equations $M$ with respect to the corresponding variables is computed.
\im Define the interaction terms   as terms in $M$ which contain variables $s\in \text{infc}$, and which, if positive, must end up in $F$.
Their complement, denoted by $V1$,  will form part of $V$.
\im A first guess  for  $F$, $F1$ is constructed as the complement of $V1$.
It contains all the interaction  terms (which involve  both disease and susceptible \com s).
\im  $F$ is obtained by  retaining  only the positive part of the matrix $F_1$, i.e. the terms which do not contain syntactic minuses.\fn[4]{we use the simplest algebraic representation of the equations, and do not study the effect which algebraic manipulations introducing minuses might have.} Finally, $V1$ is increased to $V$, which is the complement of $F$.
\im The script outputs \{M,V1,F1,F,V,K\}.
\EEN

 \begin{verbatim}
NGM[mod_,inf_]:=Module[{dyn,X,infc,M,V,F,F1,V1,K},
   dyn=mod[[1]];X=mod[[2]];
   infc=Complement[Range[Length[X]],inf];
   M=Grad[dyn[[inf]],X[[inf]]]
   (*The jacobian of the infectious equations*);
   V1=-M/.Thread[X[[infc]]->0]
   (*V1 is a first guess for V, retains all gradient terms
    which disappear when the non infectious components are null*);
   F1=M+V1
   (*F1 is a first guess for F, containing all other
   gradient terms*);
   F=Replace[F1, _. _?Negative -> 0, {2}];
   (*all terms in F1 containing minuses are set to 0*);
   V=F-M;
   K=(F . Inverse[V])/.Thread[X[[inf]]->0]//FullSimplify;
{M,V1,F1,F,V,K}]
 \end{verbatim}

  The results of this script turn out    to yield correct results in all the examples from the literature we checked in  \cite{AABJ}.

  \beR Note that the ``Replace" command in the script uses the powerful Mathematica capability of applying a ``rule" to  parts of an ``expression", specified by ``levelspec",  and that it was furnished to us by the user Michael E2 in \begin{verbatim}
https://mathematica.stackexchange.com/questions/286500/
how-to-set-to-0-all-terms-in-a-matrix-which-contain-a-minus
/287406?noredirect=1#comment715559_287406
\end{verbatim}

\eeR

\beR As an aside,   the nowadays \ME\ literature suffers a lot  from "irreproducibility", i.e. the lack of electronic notebooks to support complicated computations, which changes the  simple task of pressing enter for checking into days of work. We hope that our scripts
will help in this direction. \eeR

 \ssec{The Jacobian factorization  bound}
 Note first the following elementary fact:
 \beL A sufficient (but not necessary) condition for a polynomial with real coefficients and {\bf positive leading term}  to admit a positive
 root is that $c_0  <0$, where  $c_0 $
is the constant term   of the polynomial.

\eeL
For polynomials of degree $1$, this condition is  also necessary.
This converse result may be strenghtened  to ``Descartes polynomials".
\beD A) We will say that  a parametric polynomial with real coefficients which is irreducible over the rationals is of \Dt, if its constant coefficient may change sign, but  all the other coefficients are ``sign definite", and of the same sign (which may be supposed \wl\ to be $+$)

B) We will call a rationally reducible characteristic polynomial of \Dt, if it may be factored in  factors are either Descartes polynomials, or  linear terms  with  roots/eigenvalues which are negative (in the latter case they are irrelevant for our purpose).
\eeD

\beL A sufficient and necessary condition for a Descartes polynomial with a {\bf positive leading term}  to admit a positive
 root is that $c_0  <0$, where  $c_0 $
is the constant term   of the polynomial.

\eeL

The  Jacobian factorization approach consists in:
\BEN \im Putting all the rational factors of the characteristic polynomial of the Jacobian, in a form normalized to have positive leading term, assuming they  are sign definite (if this is not the case, this approach does not work, but may be generalized).
\im Removing all linear factors  with  eigenvalues which are negative. %, independently of the parameters.
\im For all remaining factors $F_i$ for which $c_0^{(i)} <0$ may hold for certain parameter values, rewrite this inequality into the form
$$c_0^{(i)} =c_+ - c_-=c_+(1-R_J^{(i)}) <0 \Eq R_J^{(i)}:=\fr{c_-}{c_+}>1,$$ where $c_+,c_-$ are the positive and negative parts of the expanded form of $c_0^{(i)}.$
\im Define the "Jacobian factorization $R_0$"
\begin{equation}\label{RJ} R_J=\max_i[R_J^{(i)}].\end{equation}
\EEN

%\beR We may similarly define a quantity $R_{JN}=\max_i[R_{JN}^{(i)}]$, where

\beT  A) In the instability domain, $R_J$ is a lower bound for $\inf_{F \text{ admissible}} R_F$.

 B) In the stability domain, $R_J$ is an upper bound for $\sup_{F \text{ admissible}} R_F$.

\eeT

{\Prf}   A) Fix any  admissible $F$ and let $R_F$ be its associated NGM $R_0$. Then
 $$R_J >1 \Eq \exists i: R_J^{(i)} >1 \Eq \exists i: c_0^{(i)} <0 \Lra \text{ DFE instability } \Eq R_F >1.$$  Thus
\begin{equation}\label{Rin}R_J >1  \Lra  R_F >1 \Eq R_J \leq R_F,\end{equation}
and the result follows.

B) Similar proof.

 % \red{it seems  that ${c_-}$ contains terms which are at least of degree $1$ in the $s$ variables.}

 {\bf Conjecture:} We conjecture that if all the remaining factors $F_i$ are Descartes polynomials, then $R_{J}=R_F$ for any admissible decomposition, and will denote the resulting object by $R_0$.

\beR  It holds often in simple examples that
$R_J=R_N$, or $R_J=R_N^2$, where $R_N$ denotes the $R_0$ obtained by our recipe NGM, and therefore the two $R_0$'s determine the same stability domain. %Even in the  example   in section \ref{s:Mart} where the relation between $R_J$ and $R_N$ is not that simple, they define the same stability domain.

\eeR

  %{An  example  where strict inequalities  in Theorem 1 hold} is furnished in Section \ref{s:Van12}, where the unique nonlinear  factor  of the characteristic polynomial of the Jacobian has degree $2$.
  \iffalse In that case, our Jacobian factorization algorithm may be amended.
  Indeed, the necessary and sufficient Routh-Hurwitz conditions may be written as: $c_0 <0$ or $c_0 c_1 <0$. The latter inequality, not accounted for in Thm. 1, may also be put in the form $R >1$ for some $R$, and adding that $R$ into the definition of $R_J$ improves the result.
\fi

   {\bf Open question 1:}  Under what
conditions do our NGM $R_N$ and our Jacobian  $R_J$  coincide?

We provide now the implementation  of the  Jacobian factorization approach. First, we use a utility which, for a given model, infectious set, and dummy variable (taken always as $u$, to avoid confusions), outputs the Jacobian at the DFE, the  trace
and determinant (for other purposes), the characteristic polynomial in $u$, the NGM matrix  and $R_F$.
\begin{verbatim}
JR0[mod_,inf_,u_,cn_:{}]:=
  Module[{dyn,X,par,cinf,cp,cX,jac,tr,det,chp,ngm,K,R0},
    dyn=mod[[1]];X=mod[[2]];par=mod[[3]];
    Print[" dyn=",dyn//FullSimplify//MatrixForm,X,par];
    cinf=Thread[X[[inf]]->0];
    cp=Thread[par>0];cX=Thread[X>0];
    cdfe=Join[DFE[mod,inf],cinf];
    jac=Grad[dyn,X]/.cinf/.cn;
    tr=Tr[jac];
    det=Det[jac];
    chp=CharacteristicPolynomial[jac,u];
    ngm=NGM[mod,inf];
    K=ngm[[6]];
   Print["K=",K//MatrixForm];
   R0=Assuming[Join[cp,cX],Max[Eigenvalues[K]]];
  {chp,R0,K,jac,tr,det}];
\end{verbatim}

Most of the work is done after calling this utility, by another one, JR02. This splitting of JR0 in two parts is necessary  since
the detection of the non-sign definite factors which must be analyzed  is easier to perform by eye, than to program. The JR02  script is:
\begin{verbatim}
JR02[pol_,u_]:=Module[{co,co1,cop,con,R_J},
co=CoefficientList[pol,u];
  Print["the  factor ",pol," has degree ",Length[co]-1];
  co1=Expand[co[[1]]* co[[Length[co]]]];
  Print["its leading * constant coefficient product is ",co1];
  cop=Replace[co1, _. _?Negative -> 0, {1}](*level 1 here ?*);
  con=cop-co1;
  Print["R_J is"];
  R_J=con/cop//FullSimplify;
{R_J,co}
]
\end{verbatim}

For a specific ``mod", both $R_0$'s may be obtained by typing:
\begin{verbatim}
jr = JR0[mod, inf, u];
chp = jr[[1]] // Factor
Print["factor is ", pol = chp[[5]]]
pc = JR02[pol,
   u];
   (*the script JR02 determines R_J, using the index,
    for example 5, determined by \eye inspection in
     the previous command*)
Print["R_J is ", R_J = pc[[1]] // FullSimplify]
Print["R_N is ", R_N = jr[[2]] // FullSimplify]
\end{verbatim}

To obtain the $R_0$ for our previous examples, it is enough to write down their equations and variables, for example
 \begin{verbatim}
SEIR =
    Module[{ss, ee, ii, rr, RHS, X, pars},
        ss = b - \mu s - (\beta s i+\beta_e s e) + \gamma_r r 
        - \gamma_s s;
        ee = (\beta s i+\beta_e s e) - e (\mu + e_i+e_r);
        ii = e e_i - i (\mu + \gamma + \delta);
        rr =b_r + \gamma_s s + e_r e+ \gamma i 
        - (\gamma_r + \mu) r;
        RHS = {ee, ii, rr, ss};
        X = {e, i, r, s};
        pars = Complement[Variables[RHS], X];
        {RHS, X, pars}
    ]
 \end{verbatim}
 and then call the script with the appropriate set of infectious equations ``inf".
 
  The results of this decomposition   seem to yield correct results in all the examples from the literature we checked -- see the accompanying paper \cite{AABJ}.  We conjecture, but were unable to prove,  that for \BEN \im essentially nonnegative polynomial processes having  a non-empty set of disease classes (so we deal with an epidemic); \im   processes  with a unique DFE,  at least after excluding   biologically irrelevant fixed points, like an unreachable origin;
 \im   the  local stability domain  of the DFE is non-empty, and not the full set,\EEN
   this decomposition yields ``admissible gradient decompositions", in the sense that $V^{-1}$ will contain only non-negative terms, and that it is furthermore obtainable from an equations decomposition which is admissible in the sense of \cite{Van08} (see Definition 1), and yields therefore
  the correct stability domain.

\ssec{The ``rational univariate representation" (RUR) method}\la{s:RUR}

We mention now the possibility to develop yet another  algorithm to compute a ``bifurcation $R_0$", which is based on the known fact that this parameter is expected to produce bifurcations  at $R_0=1$. Hundreds of mathematical epidemiology 
papers have already employed the idea of  reducing the fixed point system to one scalar equation in one of the disease variables, via rational substitutions for the other variables. We note here that this is a   particular case of the so-called "rational univariate representation" (RUR), but for now this is irrelevant, since RUR seems not to be implemented now in Mathematica, and we had to write our own script, included below, in which the user chooses the variable he wants to restrict to.

In the corresponding scalar polynomial the variable may be factored out (by the assumption that this is  a disease variable), and furthermore the free  coefficient of the divided polynomial may be written as $F(R_0)=G(R_0) (R_0-1)$, where  $F(R)$ is rational.

Practically, we must identify  a factor which is linear in susceptible variables, and write as a difference  of positive and negative terms. Upon  normalizing one of them to one, the other will be $R_0$, or $1/R_0$. It will  be instructive to compare the results of RUR with those of the NGM script. 
The current code for this ``budding algorithm" is
\begin{verbatim}
RUR[mod_, ind_, cn_ : {}] := Module[{dyn, X, par, eq, elim},
      dyn = mod[[1]]; X = mod[[2]]; par = mod[[3]]; 
      eq = Thread[dyn == 0]; 
      elim = Complement[Range[Length[X]],ind];
      pol = 
       Collect[GroebnerBasis[Numerator[Together[dyn /. cn]], 
         Join[par, X[[ind]]], X[[elim]]], X[[ind]]];
         ratsub = Solve[Drop[eq, ind], X[[elim]]][[1]]; 
         {ratsub, pol}
    ]
\end{verbatim}
 
\beR This approach may work only for models with demographics, in which there is a finite number of stationary fixed points,  and must fail
otherwise, for ``conservation systems" where the fixed points are only determined by the conservation of the total mass. \eeR

\sec{Endemic fixed points for some examples, via  the RUR script}
\beXa The S$I^2$R/SAIR/SEIR-FA model has a unique endemic point, with \begin{equation}\label{eeS} \s_e=\fr{\sd}{\mR_0},\mathsf e_e=(\mR_0 -1) A, A>0,\mathsf i_e=\e_e \frac{e_i}{\delta +\gamma _i+\mu }, \mathsf r_e=\rd+(\mR_0 -1) B, B>0.\end{equation}
This point is nonnegative iff $\mR_0\ge 1$, and at equality coincides with the DFE.

  \eqr{eeS} may be obtained by a symbolic solve.  In general,  for harder cases, we may attempt to  reduce the   system to one scalar polynomial equation. The following script  outputs a scalar polynomial in  a variable ``ind" specified by the user, ( chosen as $i$ for example), and also the result of solving the other variables with respect to $ind$:
\begin{verbatim}
 RUR[mod_, ind_, cn_ : {}] := Module[{dyn, X, par, eq, elim},
         dyn = mod[[1]]; X = mod[[2]]; par = mod[[3]]; 
      eq = Thread[dyn == 0];
       elim = Complement[Range[Length[X]],  ind];
         pol = 
       Collect[GroebnerBasis[Numerator[Together[dyn /. cn]], 
         Join[par, X[[ind]]], X[[elim]]], X[[ind]]];
         ratsub = Solve[Drop[eq, ind], elim][[1]]; 
         {ratsub, pol}
      ]
 \end{verbatim}

 If the RHS of the infectious equation has not been simplified by $i$, the polynomial must further be divided by it. Finally, this yields a polynomial of degree $1$, and
 \begin{verbatim}
 in={1};pol=Grobpol[mode,in];cof = CoefficientList[pol, i];
 cof[[1]];
 po=Sum[cof[[k]] e^{k-2},{k,Length[cof]}];
 Solve[po==0,e]//FullSimplify
 \end{verbatim}
 will reveal the complicated formula of $A$ in \eqr{eeS}.
 \eeXa

\sec{The universal language of  pseudo-linear, essentially non-negative, {mass action} dynamical systems}\la{s:back}

\ssec{Pseudo-linearity}
The  sister disciplines of   mathematical epidemiology, chemical reaction networks (CRN), ecology, virology, biochemical systems, etc., started all as  collections of examples of
   ``pseudo-linear"
differential systems {\bf parametrized by two matrices} $\Sb , {\Yb}$:
\beq \la {Wdef} && \dot{\xb}=f(\xb)=\sum_{k=1}^{n_R} {\mathbf{s}}_k \xb^{\yb_k}=\Sb \xb^{\Yb}, \quad  \xb, \yb_k, {\mathbf{s}}_k \in \R_+^{n \times 1}, \eeq
where $\xb^{\Yb} \in \R_+^{n_R \times 1}$ is a column vector of monomials, $\Yb \in \R^{n \times n_R}$ is the ``matrix of source exponents"
 and $\Sb\in \R^{n \times n_R}$ is the ``stoichiometric matrix of  direction vectors" (formed \resp\ by joining exponents $\yb_1,...,\yb_{n_R}$ and directions ${\mathbf{s}}_1,...,{\mathbf{s}}_{n_R}$ as columns). Note that any polynomial dynamical systems can be uniquely written in such form, for some distinct $\yb_i$ , and
 non-zero ${\mathbf{s}}_i$ \cite{craciun2022algorithm},  that $\Yb$ is not unique, and that its  dimension may be easily increased .

\ssec{Essential non-negativity and mass action representation}

 Kinetic systems must be ``essentially non-negative",   meaning that they leave invariant the nonnegative orthant.%\fn[3]{A more restrictive possible definition of kinetic systems could be that of  ``essentially polytopal  ODE systems", i.e. systems which may not exit a bounded polytope.}

\beR An obvious sufficient condition for the essential non-negativity (i.e. the preservation of the
nonnegative octant) of a polynomial system $X'=f(X)$ is that each component $f_i(X)$ may be decomposed
   as \begin{equation}\label{Hc} f_i(X)=g_i(X)-x_i h_i(X),\end{equation}
  where $g_i,  h_i$ are polynomials with nonnegative coefficients, i.e. if all negative terms in an equation  contain the variable whose rate is  given by the equation.

   \beD Terms which do not satisfy \eqr{Hc} are called negative cross-effects. \eeD

\eeR

  \beXa The  Lorentz system (a famous example of chaotic behavior)
$$\bc x'=\sigma(y-x)\\
y'=\rho x- y-x z\\z'=xy -\beta z\ec$$
 does not satisfy \eqr{Hc}, due to the $-x z$ term in the $y$ equation.
 \eeXa

 The following result, sometimes called the ``Hungarian lemma" is well-known in  the \CRN\ literature  \cite{hun},  \cite[Thm 6.27]{Toth}:

 \beL  \la{l:hun} A polynomial system admits an essentially non-negative ``mass-action" representation iff
 there are no negative cross-effects, i.e. if \eqr{Hc} holds.
\eeL

\iffalse
\beR
One characterization of  mass action chemical reaction networks is as the subclass of pseudo-linear systems  which
admit a representation
\be{pr} \Sb=\Yb \Ab,\ee
where $\Ab$ is a Laplacian matrix. There exists a canonical representation
 which achieves this.  Unfortunately, it yields a  $\Yb$ with large dimension, and several following works have  produced minimal representations in various senses. Whether these representations are small enough to be useful given the current limitations
of our computers is not yet clear.
\eeR
\fi

\section{Examples of epidemic models presented as chemical reaction networks systems \la{s:3}}
\ssec{A brief perspective of chemical reaction networks theory}\la{s:CRN}
CRN researchers succeeded in developing several general laws under assumptions like {\bf ''low deficiency"} and  {\bf '' complex-balancedness/toricity"} -- see the pioneering works of Feinberg, Horn -Jackson and Clarke\cite{clarke1993method}, continued by \cite{craciun2005multiple,Gater,mincheva2007graph,mincheva2008multigraph,craciun2008homotopy,
CraStu,mincheva2012oscillations,kaltenbach2012unified,otero2012characterizing,
otero2017chemical,
craciun2020delay,yu2022graph,domijan2009bistability,errami2012computing,
errami2015detection,england2017symbolic,Licht,
sadeghimanesh2022resultant}, among others. These results are often associated to   specific graph structures associated to a ``Laplacian matrix". %There do not seem to exist for now applications of these results outside the field.

A significant recent development was the discovery by Karin Gatermann -- see \fe\ \cite{Gater}--of the connection between the graphs of mass-action kinetics and algebra. Gatermann, and subsequently   Craciun,  Dickenstein,  Shiu, and collaborators  showed that graph  concepts like detailed balance (one of the first concerns of \CRNs, since Boltzmann), may be translated in algebraic language into the binomiality of the ideal generated by the fixed point equations.
Subsequently, \cite{Gater} observed  that the
Jacobian  may be viewed as  a weighted adjacency matrix of a directed graph, whose cycles had appeared already in the necessary ``sign-stability" conditions in the Quirk–Ruppert–Maybee theorem \cite{jeffries1977matrix}.

  %One of the main points of their analysis was the change to ``reaction variables", after which the ``endemic point" with positive coordinates becomes unique \cite{Gater}.

%\iffalse

 \beR  {\bf One remarkable result} of CRN theory is that under a certain relation between
 the parameters which used to be called Wegscheider condition and is now called toricity,
 a mass-action system is ''{\bf complex-balanced/toric}", i.e. there exists a unique fixed point \satg\
 $\Ab \xb^{\Yb}=0$),  within each fiber
 determined by the constraints.
 \eeR

 Chemical reaction networks theory  contains both numerous  general results, and also an impressive number of open problems. One examples
 of the latter is:

 {\bf Q}:  Can chaos occur in three dimensional ``bimolecular" (with linear and quadratic interactions only) essentially nonnegative mass action models? All the examples we are aware of are either four-dimensional, or non-polynomial-- see for example \cite{letellier2013can}, or non-quadratic--
 see \cite[Table 1]{Wilh} for a review of some  minimal mass action systems admitting chaos.
 \iffalse
 {\bf Q1:} Can  ``complex balance/toricity" be defined
for arbitrary kinetic systems  directly in terms of the matrices $\Sb,\Yb$ in the definition, without constructing first associated  graphs (which are not unique),  bypassing thus the Laplacian decomposition which establishes that a kinetic system  is a  \CRN?
\fi

%\input{CRN}
%\sec{Appendix}
\subsection{Can chemical reaction network  (CRN) theory  help for solving epidemic models?} \la{s:SIRCRN}

The answer is yes, as shown recently for Capasso type SIR models and SIRnS models, in \cite{VAA}. Here, we pass briefly in review
some other CRN tools which might turn out useful in the future. 

 To wet the appetite for this question, we explain now    how the seven parameter kinetic system SIR/V+S  from \cite{AABBGH} may be represented as a CRN.
 With $ X= \bep \s \\ \mathsf i \\ \mathsf r\eep$, this
  may be written    as

  \begin{align}
& X'(t)= \Lambda \mathsf i \bep 1 \\-1\\0 \eep+ (\b-\delta) \s  \mathsf i \bep -1\\1\\0\eep +(\beta_r-\delta) \mathsf i   \mathsf r \bep 0\\1\\-1\eep
\\&+ \g  \mathsf i  \bep 0\\-1\\1\eep+ \g_s \s \bep -1\\0\\1\eep
+\g_r \mathsf r \bep 1\\0\\-1\eep +\Lambda \mathsf r \bep 1\\0\\-1\eep.% + \delta \s   \mathsf i  \bep 1\\0\\0\eep+ \delta  \i^2  \bep 0\\1\\0\eep.
\no
\end{align}

Assuming now $\delta<\beta, \delta<\beta_r$, the positivity of the rate constants allows us to introduce a ``stoichiometry matrix"
$$\Sb= \bep 1 \\-1\\0 \eep   \bep -1\\1\\0\eep \bep 0\\1\\-1\eep \bep 0\\-1\\1\eep
\bep -1\\0\\1\eep \bep 1 \\0\\-1 \eep \bep 1 \\0\\-1 \eep .   $$

\beR The stoichiometry matrix has  rank 2, same as $ \bep 1 \\-1\\0 \eep
\bep 0\\-1\\1\eep \bep -1\\0\\1\eep,   $
 since $\text{row}(1)+\text{row}(2)=- \text{row}(3) $.
   %\bep 1\\0\\0\eep  \bep 0\\1\\0\eep  \bep 0\\0\\1\eep

\eeR

The reactions of the corresponding CRN, assuming ``mass-action form", may be written as $$\bc I->S & \Lambda \mathsf i \\I->R &\g  \mathsf i\\ S+I ->2 I &(\beta-\delta) \s \mathsf i\quad (\mbox{with}\; \beta>\delta)\\I+R-> 2I &(\beta_r-\delta) \mathsf r \mathsf i \quad (\mbox{with}\; \beta_r>\delta)\\ S->R &\g_s \s\\  R->S &\g_r  \mathsf r\\ R->S &\Lambda \; \mathsf r\ec.$$
% S + I->2S+I &(\delta \s \i)\\2 I->3I &(\delta  \i^2)\\.$$
It has  $n_V=6$ vertices $S,I,R,S+I,S+R,2 I$ (the first $n_s=5$ of which are sources), $n_R=7$  edges,
and 2 linkage classes.  The deficiency of the CRN is $\delta=n_V-rank(S)-n_C=6-2-2=2$. It may be possible to  produce a representation of lower deficiency,  and thus take advantage of classic CRN results.

\figu{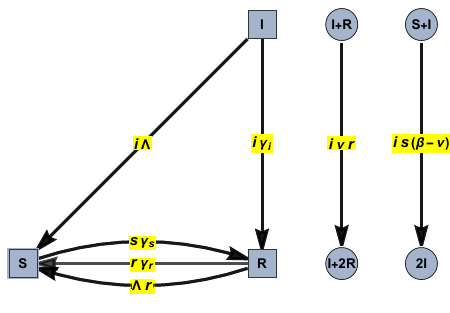}{The {Feinberg-Horn-Jackson graph} of the ``SIR network" with $7$ reactions.}{1}

\beXa Let us consider now the simpler {five} parameters  SIR/V+S model obtained by taking $\beta_r=0=\delta$.%, and $\delta=0$ in the quadratic terms.
The reactions of the corresponding mass-action CRN are
 $$\bc I->S & \Lambda \mathsf i \\I->R &\g  \mathsf i\\ S+I ->2 I &\beta \s \mathsf i
 \\ S->R &\g_s \s\\  R->S &\g_r  \mathsf r\\ R->S &\Lambda \; \mathsf r\ec.$$

It has  $n_V=5$ vertices , $n_R=6$  edges,
and 2 linkage classes.  The deficiency of the CRN is $\delta=n_V-rank(S)-n_C=5-2-2=1$.

\eeXa

\iffalse 
\documentclass{article}
%%%%%%%%%%%%%%%%%%%%%%%%%%%%%%%%%%%%%%%%%%%%%%%%%%%%%%%%%%%%%%%%%%%%%%%%%%%%%%%%%%%%%%%%%%%%%%%%%%%%%%%%%%%%%%%%%%%%%%%%%%%%%%%%%%%%%%%%%%%%%%%%%%%%%%%%%%%%%%%%%%%%%%%%%%%%%%%%%%%%%%%%%%%%%%%%%%%%%%%%%%%%%%%%%%%%%%%%%%%%%%%%%%%%%%%%%%%%%%%%%%%%%%%%%%%%
\usepackage{amsfonts}

%TCIDATA{OutputFilter=Latex.dll}
%TCIDATA{Version=5.00.0.2552}
%TCIDATA{<META NAME="SaveForMode" CONTENT="1">}
%TCIDATA{LastRevised=Tuesday, June 13, 2023 18:27:42}
%TCIDATA{<META NAME="GraphicsSave" CONTENT="32">}

\input{tcilatex}

\begin{document}
\fi 

\sec{Lagrange-Hamilton geometric objects for general dynamical systems}
\la{s:Neag}

For the manifold $M=\mathbb{R}^{n+1}$, whose coordinates are%
\[
\left( x^{k}\right) _{k=\overline{1,n+1}}=\left( \mathbf{i},s\right) ,
\]%
then we denote by $TM$ (respectively $T^{\ast }M$) the tangent (respectively
cotangent) bundle, whose coordinates are $\left( x^{i},y^{i}\right) _{i=%
\overline{1,n+1}}$ (respectively $\left( x^{i},p_{i}\right) _{i=\overline{%
1,n+1}}).$ 

The solutions of class $C^{2}$ of the dynamical system (\ref{DS}) are the
global minimum points for the \textit{least squares Lagrangian} $%
L:TM\rightarrow \mathbb{R}$ (see \cite{Udr}, \cite{Nea-Udr}) given by%
\begin{equation}
L(x,y)=\left( y^{1}-X^{1}(x)\right) ^{2}+\left( y^{2}-X^{2}(x)\right)
^{2}+...+\left( y^{n+1}-X^{n+1}(x)\right) ^{2}.  \label{LSL}
\end{equation}%
Consequently, via the Euler-Lagrange equations\footnote{%
The Einstein convention of summation is adopted all over this Section.}%
\[
\frac{\partial L}{\partial x^{k}}-\frac{d}{dt}\left( \frac{\partial L}{
\partial y^{k}}\right) =0,\quad y^{k}=\frac{dx^{k}}{dt},\quad k=\overline{
1,n+1},
\]%
\begin{equation}
\Leftrightarrow \frac{d^{2}x^{k}}{dt^{2}}+2G^{k}(x,y)=0,\quad G^{k}=\frac{1}{
2}\left( \frac{\partial ^{2}L}{\partial x^{j}\partial y^{k}}y^{j}-\frac{
\partial L}{\partial x^{k}}\right) ,  \label{E-L}
\end{equation}%
we can build a whole natural collection of nonzero Lagrangian geometrical
objects (such as nonlinear connection, d-torsions and Yang-Mills
electromagnetic-like energy) that characterize the dynamical system (\ref{DS}%
). A detailed exposition of the Lagrange geometry on tangent bundles,
together with the Lagrangian least squares variational method for dynamical
systems, find in the works: Miron and Anastasiei \cite{Mir-An}, Udri\c{s}te
and Neagu \cite{Udr}, \cite{Nea-Udr}, and Balan and Neagu \cite{Bal-Nea}.

Also, we can construct the \textit{least squares Hamiltonian }$H:T^{\ast
}M\rightarrow \mathbb{R}$, associated with the Lagrangian (\ref{LSL}), which
is defined by%
\begin{equation}
H(x,p)=\displaystyle{\frac{1}{4}}\left(
p_{1}^{2}+p_{2}^{2}+...+p_{n+1}^{2}\right)
+X^{1}(x)p_{1}+X^{2}(x)p_{2}+...+X^{n+1}(x)p_{n+1},  \label{LSH}
\end{equation}%
where $p_{r}=\partial L/\partial y^{r}$ and $H=p_{r}y^{r}-L$. It follows
that we can construct again a natural and distinct collection of nonzero
Hamiltonian geometrical objects (such as nonlinear connection and
d-torsions), which also characterize the dynamical system (\ref{DS}). The
detailed construction of the Hamilton geometry on cotangent bundles,
together with the Hamiltonian least squares variational method for dynamical
systems, can be found in the monographs: Miron et al. \cite%
{Miron-Hr-Shim-Sab} and Neagu and Oan\u{a} \cite{Nea-Oana}.

The Lagrange-Hamilton geometries produced by the Lagrangian (\ref{LSL}) and
Hamiltonian (\ref{LSH}) are achieved via the nonzero geometrical objects
(see the works \cite{Bal-Nea}, \cite{Nea-Oana} and \cite{Nea-Lit}):

\begin{enumerate}
\item $\mathcal{N}=\left( N_{j}^{i}\right) _{i,j=\overline{1,n+1}}=- %
\displaystyle{\frac{1}{2}}\left[ J-J^{t}\right] $ -- is the \textit{\
Lagrangian nonlinear connection} on the tangent bundle $TM$, where $%
N_{j}^{i}=\partial G^{i}/\partial y^{j};$\medskip

\item $R_{k}=\left( R_{jk}^{i}\right) _{i,j=\overline{1,n+1}}=%
\displaystyle 
\frac{\partial \mathcal{N}}{\partial x^{k}},$ $\forall $ $k=\overline{1,n+1}
, $ -- are the \textit{Lagrangian d-torsions}, where 
\[
R_{jk}^{i}={\frac{\delta N_{j}^{i}}{\delta x^{k}}}-\frac{\delta N_{k}^{i}}{
\delta x^{j}},\quad {\frac{\delta }{\delta x^{k}}={\frac{\partial }{\partial
x^{k}}}-N_{k}^{r}{\frac{\partial }{\partial y^{r}};}} 
\]

\item $\mathcal{EYM}(x)=\displaystyle{\frac{1}{2}}\cdot Trace\left[ F\cdot
F^{t}\right] $, where $F=-\mathcal{N}$, -- is the \textit{Lagrangian
Yang-Mills e\-lec\-tro\-mag\-ne\-tic-like energy};\medskip

\item $\mathbf{N}=\left( N_{ij}\right) _{i,j=\overline{1,n+1}}=J+J^{t}$ --
is the \textit{Hamiltonian nonlinear connection} on the cotangent bundle $%
T^{\ast }M$ expressed by 
\[
N_{ij}=\frac{\partial ^{2}H}{\partial x^{j}\partial p_{i}}+\frac{\partial
^{2}H}{\partial x^{i}\partial p_{j}}; 
\]

\item $\mathbf{R}_{k}=\left( R_{kij}\right) _{i,j=\overline{1,n+1}}= %
\displaystyle\frac{\partial }{\partial x^{k}}\left[ J-J^{t}\right] =-2R_{k},$
$\forall $ $k=\overline{1,n+1},$ -- are the \textit{Hamiltonian d-torsions }
defined by 
\[
R_{kij}={\frac{\delta N_{ki}}{\delta x^{j}}}-{\frac{\delta N_{kj}}{\delta
x^{i}}},\quad {\frac{\delta }{\delta x^{j}}={\frac{\partial }{\partial x^{j}}
}-N_{rj}{\frac{\partial }{\partial p_{r}},}} 
\]
where $J=\left( \partial X^{i}/\partial x^{j}\right) _{i,j=\overline{1,n+1}}$
is the Jacobian matrix of $X.$
\end{enumerate}

\noindent \textbf{Open problem.} Are there applications in mathematical
epidemiology or other disciplines dealing with kinetic systems for the above
Lagrange-Hamilton geometrical objects naturally associated to any dynamical
system?

%\end{document}

\ssec{Lagrange-Hamilton geometric objects for SIR-PH-FA   epidemic models}\la{s:n1}

We  compute here the  geometrical objects from Section \ref{s:Neag}, starting by the particular
case of SIR without loss of immunity, which is  a SI-PH-FA model with  $n=1$,  $V=b+\gamma $ and $B=\beta $ (see Example \ref{ex:SEIR}). It follows
that we have the coordinates%
\[
\left( x^{1},x^{2}\right) =\left( \mathtt{i,s}\right)
\]%
on $\mathbb{R}^{2}$, where
\[
\mathtt{i}(t)=\frac{I(t)}{N(t)}=\frac{I(t)}{I(t)+S(t)+R(t)},\quad \mathtt{s}%
(t)=\frac{S(t)}{N(t)}=\frac{S(t)}{I(t)+S(t)+R(t)},
\]%
together with%
\[
X^{1}\left( \mathtt{i,s}\right) =\beta \mathtt{si-}\left( b+\gamma \right)
\mathtt{i,\quad }X^{2}\left( \mathtt{i,s}\right) =b-\left( b+\gamma _{%
\mathtt{s}}\right) \mathtt{s}-\beta \mathtt{si,}
\]%
where: (i) $\beta $ is the transmission rate by which the infectious infect
susceptible people; (ii) $1/\gamma $ is the average time an infectious
individual may infect others; (iii) $b$ is the birth rate on capita on a
unit of time. Consequently, the corresponding Jacobian is given by%
\[
\mathtt{J}=\left(
\begin{array}{cc}
\mathtt{s}\beta -\left( b+\gamma \right) & \beta \mathtt{i} \\
-\mathtt{s}\beta & -(b+\gamma _{\mathtt{s}})-\beta \mathtt{i}%
\end{array}%
\right) .
\]%
In conclusions, by simple computations, we obtain the following
Lagrange-Hamilton geometrical objects:

\begin{enumerate}
\item the Lagrangian nonlinear connection skew-symmetric matrix on the
tangent bundle $T\mathbb{R}^{2}$:
\[
\mathtt{N}=-\displaystyle{\frac{1}{2}}\left[ \mathtt{J}-\mathtt{J}^{t}\right]
=-\frac{1}{2}\left(
\begin{array}{cc}
0 & \beta \mathtt{i}+\beta \mathtt{s} \\
-\beta \mathtt{i}-\beta \mathtt{s} & 0%
\end{array}
\right) ;
\]

\item The Lagrangian d-torsion skew-symmetric matrices:
\[
\mathtt{R}_{1}=\frac{\partial \mathtt{N}}{\partial \mathtt{i}}=-\frac{1}{2}
\left(
\begin{array}{cc}
0 & \beta \\
-\beta & 0%
\end{array}
\right) ;
\]
\[
\mathtt{R}_{2}=\frac{\partial \mathtt{N}}{\partial \mathtt{s}}=-\frac{1}{2}
\left(
\begin{array}{cc}
0 & \beta \\
-\beta & 0%
\end{array}
\right) ;
\]

\item the Lagrangian Yang-Mills electromagnetic-like energy:
\[
\mathcal{EYM}(\mathtt{i},\mathtt{s})=\frac{\beta ^{2}}{4}\left( \mathtt{i}+%
\mathtt{s}\right) ^{2};
\]

\item the Hamiltonian nonlinear connection symmetric matrix on the cotangent
bundle $T^{\ast }\mathbb{R}^{2}$:
\[
\mathsf{N}=\mathtt{J}+\mathtt{J}^{t}=\left(
\begin{array}{cc}
2\beta \mathtt{s}-2\left( b+\gamma \right) & \beta \mathtt{i}-\beta \mathtt{s%
} \\
\beta \mathtt{i}-\beta \mathtt{s} & -2(b+\gamma _{\mathtt{s}})-2\beta
\mathtt{i}%
\end{array}%
\right) ;
\]

\item the Hamiltonian d-torsion matrices are $\boldsymbol{R}_{k}=-2\mathtt{R}%
_{k},$ $\forall $ $k=\overline{1,2}.$
\end{enumerate}

We compute now the preceding Lagrange-Hamilton geometrical objects for the
general class of SIR-PH-FA epidemic models (for proportions), introduced in
[AAB+22]. Let us write the SIR-PH-FA model with $\gamma _{r}=0$  as a dynamical system%
\begin{equation}
\frac{dx^{i}}{dt}=X^{i}(x(t)),\quad i=\overline{1,n+1},  \label{DS}
\end{equation}%
where the vector field $X=\left( X^{i}(x)\right) _{i=\overline{1,n+1}}$ on $%
M $ is given in the SIR-PH-FA equation   \eqr{SYRPH} above. Recall that the Jacobian for SIR-PH-FA model without loss of
immunity ($\gamma _{r}=0$),  given in \eqr{JacPh}, is:

\begin{equation}
J=\left(
\begin{array}{cc}
sB^{t}-V^{t} & B^{t}i^{t} \\
-s\beta^{t} & -(b+\gamma _{\mathtt{s}})-\beta^{t}\mathbf{%
i}^{t}%
\end{array}%
\right) ,  \label{Jacobian_SIR-PH-FA}
\end{equation}%
where $i=\left( i_{1},i_{2},...,i_{n}\right) $ and $\beta%
=\left(
\begin{array}{c}
\beta _{1} \\
\beta _{2} \\
\vdots \\
\beta _{n}%
\end{array}%
\right) $. Consequently, the corresponding formulas for the general class of
SIR-PH-FA epidemic models, via the coordinates%
\[
\left( x^{1},x^{2}\right) =\left( i\mathtt{,s}\right)
\]%
on $\mathbb{R}^{n+1}$ are obtained simply by replacing $b+\gamma $ by the
matrix $V$, and $\beta $ by the column $\beta$:

\begin{enumerate}
\item the Lagrangian nonlinear connection $\mathfrak{N}=\left( \mathfrak{N}%
_{j}^{i}\right) _{i,j=\overline{1,n+1}}$ on $T\mathbb{R}^{n+1}$:
\[
\mathfrak{N}=-\displaystyle{\frac{1}{2}}\left[ J-J^{t}\right] =-\frac{1}{2}%
\left(
\begin{array}{cc}
\mathtt{s}(B^{t}-B)-(V^{t}-V) & B^{t}i^{t}+\mathtt{s}\beta
\\
-\mathtt{s}\beta^{t}-iB & 0%
\end{array}%
\right) ;
\]

\item The Lagrangian d-torsion matrices:
\[
\mathfrak{R}_{k}=\frac{\partial \mathfrak{N}}{\partial i_{k}}=-\frac{1}{2}%
\left(
\begin{array}{cc}
\mathbb{O} & B_{k}^{t} \\
-B_{k} & 0%
\end{array}%
\right) ,\quad k=\overline{1,n};
\]%
\[
\mathfrak{R}_{n+1}=\frac{\partial \mathfrak{N}}{\partial \mathtt{s}}=-\frac{1%
}{2}\left(
\begin{array}{cc}
B^{t}-B & \beta \\
-\beta^{t} & 0%
\end{array}%
\right) ,
\]%
where $B_{k}=\left( B_{k,1},B_{k,2},...,B_{k,n}\right) $;

\item the Lagrangian Yang-Mills electromagnetic-like energy:
\[
\mathcal{EYM}(i,\mathtt{s})=\sum_{r=1}^{n}\sum_{q=r+1}^{n+1}\left(
\mathfrak{N}_{r}^{q}\right) ^{2};
\]

\item the Hamiltonian nonlinear connection $T^{\ast }\mathbb{R}^{n+1}$:
\[
\mathbb{N}=J+J^{t}=\left(
\begin{array}{cc}
\mathtt{s}(B^{t}+B)-(V^{t}+V) & B^{t}i^{t}-\mathtt{s}\beta
\\
iB-\mathtt{s}\beta^{t} & -2(b+\gamma _{\mathtt{s}})-2%
i\beta %
\end{array}%
\right) ;
\]

\item the Hamiltonian d-torsion matrices are $\mathbb{R}_{p}=-2\mathfrak{R}%
_{p},$ $\forall $ $p=\overline{1,n+1}.$
\end{enumerate}

 {\bf Acknowledgement} We thank Dan Goreac, Daniel Lichtblau, %Florian Nill
 and Janos Toth for useful
 advice.
%\bibliographystyle{unsrt}
  %\bibliographystyle{amsalpha}
 % \section*{References}
\bibliographystyle{amsalpha}
\bibliography{Pare40}

\newcommand{\etalchar}[1]{$^{#1}$}
\providecommand{\bysame}{\leavevmode\hbox to3em{\hrulefill}\thinspace}
\providecommand{\MR}{\relax\ifhmode\unskip\space\fi MR }
% \MRhref is called by the amsart/book/proc definition of \MR.
\providecommand{\MRhref}[2]{%
  \href{http://www.ams.org/mathscinet-getitem?mr=#1}{#2}
}
\providecommand{\href}[2]{#2}
\begin{thebibliography}{ABvdD{\etalchar{+}}07}

\bibitem[AAB{\etalchar{+}}23]{AABBGH}
Florin Avram, Rim Adenane, Lasko Basnarkov, Gianluca Bianchin, Dan Goreac, and
  Andrei Halanay, \emph{An age of infection kernel, an {R} formula, and further
  results for arino--brauer {A}, {B} matrix epidemic models with varying
  populations, waning immunity, and disease and vaccination fatalities},
  Mathematics \textbf{11} (2023), no.~6, 1307.

\bibitem[AABH22]{AABH}
Florin Avram, Rim Adenane, Gianluca Bianchin, and Andrei Halanay,
  \emph{Stability analysis of an {E}ight parameter {SIR}-type model including
  loss of immunity, and disease and vaccination fatalities}, Mathematics
  \textbf{10} (2022), no.~3.

\bibitem[AABJ23]{AABJ}
Florin Avram, Rim Adenane, Lasko Basnarkov, and Matthew~D Johnston,
  \emph{Algorithmic approach for a unique definition of the next-generation
  matrix}, Mathematics \textbf{12} (2023), no.~1, 27.

\bibitem[AAGH23]{AAGH}
Florin Avram, Rim Adenane, Dan Goreac, and Andrei Halanay, \emph{Explicit
  mathematical epidemiology results on age renewal kernels and r0 formulas are
  often consequences of the rank one property of the next generation matrix},
  arXiv preprint arXiv:2307.04774 (2023).

\bibitem[AAGH24]{AABGH}
\bysame, \emph{Does mathematical epidemiology have general laws, besides the
  {DFE} stability theorem?}, Monografías Matematicas Garcia de Galdeano
  \textbf{43} (2024), 11--20.

\bibitem[AAH22]{AAH}
Florin Avram, Rim Adenane, and Andrei Halanay, \emph{New results and open
  questions for sir-ph epidemic models with linear birth rate, waning immunity,
  vaccination, and disease and vaccination fatalities}, Symmetry \textbf{14}
  (2022), no.~5, 995.

\bibitem[ABG20]{Ando}
Alessia And{\`o}, Dimitri Breda, and Giulia Gava, \emph{How fast is the linear
  chain trick? a rigorous analysis in the context of behavioral epidemiology.},
  Mathematical Biosciences and Engineering \textbf{17} (2020), no.~5,
  5059--5085.

\bibitem[ABvdD{\etalchar{+}}07]{Arino}
Julien Arino, Fred Brauer, Pauline van~den Driessche, James Watmough, and
  Jianhong Wu, \emph{A final size relation for epidemic models}, Mathematical
  Biosciences \& Engineering \textbf{4} (2007), no.~2, 159.

\bibitem[AKK{\etalchar{+}}20]{Ansumali}
Santosh Ansumali, Shaurya Kaushal, Aloke Kumar, Meher~K Prakash, and
  M~Vidyasagar, \emph{Modelling a pandemic with asymptomatic patients, impact
  of lockdown and herd immunity, with applications to sars-cov-2}, Annual
  reviews in control (2020).

\bibitem[AP20]{arino2020simple}
Julien Arino and St{\'e}phanie Portet, \emph{A simple model for covid-19},
  Infectious Disease Modelling \textbf{5} (2020), 309--315.

\bibitem[BDDG{\etalchar{+}}12]{Breda}
Dimitri Breda, Odo Diekmann, WF~De~Graaf, A~Pugliese, and R~Vermiglio, \emph{On
  the formulation of epidemic models (an appraisal of kermack and mckendrick)},
  Journal of biological dynamics \textbf{6} (2012), no.~sup2, 103--117.

\bibitem[BELE19]{brouwer2019phenotypic}
Andrew~F Brouwer, Marisa~C Eisenberg, Nancy~G Love, and Joseph~NS Eisenberg,
  \emph{Phenotypic variations in persistence and infectivity between and within
  environmentally transmitted pathogen populations impact population-level
  epidemic dynamics}, BMC Infectious Diseases \textbf{19} (2019), 1--13.

\bibitem[BN11]{Bal-Nea}
Vladimir Balan and Mircea Neagu, \emph{Jet single-time lagrange geometry and
  its applications}, John Wiley \& Sons, 2011.

\bibitem[Bra05]{Brauer05}
Fred Brauer, \emph{The kermack--mckendrick epidemic model revisited},
  Mathematical biosciences \textbf{198} (2005), no.~2, 119--131.

\bibitem[CCH18]{cassidy2018recipe}
Tyler Cassidy, Morgan Craig, and Antony~R Humphries, \emph{A recipe for state
  dependent distributed delay differential equations}, arXiv preprint
  arXiv:1811.05930 (2018).

\bibitem[CDE18]{Champredon}
David Champredon, Jonathan Dushoff, and David~JD Earn, \emph{Equivalence of the
  erlang-distributed seir epidemic model and the renewal equation}, SIAM
  Journal on Applied Mathematics \textbf{78} (2018), no.~6, 3258--3278.

\bibitem[CDSS09]{CraStu}
Gheorghe Craciun, Alicia Dickenstein, Anne Shiu, and Bernd Sturmfels,
  \emph{Toric dynamical systems}, Journal of Symbolic Computation \textbf{44}
  (2009), no.~11, 1551--1565.

\bibitem[CF05]{craciun2005multiple}
Gheorghe Craciun and Martin Feinberg, \emph{Multiple equilibria in complex
  chemical reaction networks: I. the injectivity property}, SIAM Journal on
  Applied Mathematics \textbf{65} (2005), no.~5, 1526--1546.

\bibitem[CHW08]{craciun2008homotopy}
Gheorghe Craciun, J~William Helton, and Ruth~J Williams, \emph{Homotopy methods
  for counting reaction network equilibria}, Mathematical biosciences
  \textbf{216} (2008), no.~2, 140--149.

\bibitem[CJ93]{clarke1993method}
Bruce~L Clarke and Weimin Jiang, \emph{Method for deriving hopf and saddle-node
  bifurcation hypersurfaces and application to a model of the
  belousov--zhabotinskii system}, The Journal of chemical physics \textbf{99}
  (1993), no.~6, 4464--4478.

\bibitem[CJY22]{craciun2022algorithm}
Gheorghe Craciun, Jiaxin Jin, and Polly~Y Yu, \emph{An algorithm for finding
  weakly reversible deficiency zero realizations of polynomial dynamical
  systems}, arXiv preprint arXiv:2205.14267 (2022).

\bibitem[CMPP20]{craciun2020delay}
Gheorghe Craciun, Maya Mincheva, Casian Pantea, and Y~Yu Polly, \emph{Delay
  stability of reaction systems}, Mathematical Biosciences \textbf{326} (2020),
  108387.

\bibitem[DGM18]{Diek18}
Odo Diekmann, Mats Gyllenberg, and JAJ Metz, \emph{Finite dimensional state
  representation of linear and nonlinear delay systems}, Journal of Dynamics
  and Differential Equations \textbf{30} (2018), no.~4, 1439--1467.

\bibitem[DHB13]{DiekHeesBrit}
O.~Diekmann, H.~Heesterbeek, and T.~Britton, \emph{{Mathematical Tools for
  Understanding Infectious Disease Dynamics}}, Princeton Univ. Press, 2013.

\bibitem[DHM90]{Diek}
Odo Diekmann, Johan Andre~Peter Heesterbeek, and Johan~AJ Metz, \emph{On the
  definition and the computation of the basic reproduction ratio r0 in models
  for infectious diseases in heterogeneous populations}, Journal of
  mathematical biology \textbf{28} (1990), no.~4, 365--382.

\bibitem[DHR10]{Diek10}
Odo Diekmann, JAP Heesterbeek, and Michael~G Roberts, \emph{The construction of
  next-generation matrices for compartmental epidemic models}, Journal of the
  Royal Society Interface \textbf{7} (2010), no.~47, 873--885.

\bibitem[DI22]{Diek22}
Odo Diekmann and Hisashi Inaba, \emph{A systematic procedure for incorporating
  separable static heterogeneity into compartmental epidemic models}, arXiv
  preprint arXiv:2207.02339 (2022).

\bibitem[DK09]{domijan2009bistability}
Mirela Domijan and Markus Kirkilionis, \emph{Bistability and oscillations in
  chemical reaction networks}, Journal of Mathematical Biology \textbf{59}
  (2009), no.~4, 467--501.

\bibitem[EEG{\etalchar{+}}15]{errami2015detection}
Hassan Errami, Markus Eiswirth, Dima Grigoriev, Werner~M Seiler, Thomas Sturm,
  and Andreas Weber, \emph{Detection of hopf bifurcations in chemical reaction
  networks using convex coordinates}, Journal of Computational Physics
  \textbf{291} (2015), 279--302.

\bibitem[EEG{\etalchar{+}}17]{england2017symbolic}
Matthew England, Hassan Errami, Dima Grigoriev, Ovidiu Radulescu, Thomas Sturm,
  and Andreas Weber, \emph{Symbolic versus numerical computation and
  visualization of parameter regions for multistationarity of biological
  networks}, International Workshop on Computer Algebra in Scientific
  Computing, Springer, 2017, pp.~93--108.

\bibitem[ESEW12]{errami2012computing}
Hassan Errami, Werner~M Seiler, Markus Eiswirth, and Andreas Weber,
  \emph{Computing hopf bifurcations in chemical reaction networks using
  reaction coordinates}, International Workshop on Computer Algebra in
  Scientific Computing, Springer, 2012, pp.~84--97.

\bibitem[Fen07]{Feng}
Zhilan Feng, \emph{Final and peak epidemic sizes for {SEIR} models with
  quarantine and isolation}, Mathematical Biosciences \& Engineering \textbf{4}
  (2007), no.~4, 675.

\bibitem[GES05]{Gater}
Karin Gatermann, Markus Eiswirth, and Anke Sensse, \emph{Toric ideals and graph
  theory to analyze hopf bifurcations in mass action systems}, Journal of
  Symbolic Computation \textbf{40} (2005), no.~6, 1361--1382.

\bibitem[Gun03]{gunawardena2003chemical}
Jeremy Gunawardena, \emph{Chemical reaction network theory for in-silico
  biologists}, Notes available for download at http://vcp. med. harvard.
  edu/papers/crnt. pdf \textbf{5} (2003).

\bibitem[HCH10]{Haddad}
Wassim~M Haddad, VijaySekhar Chellaboina, and Qing Hui, \emph{Nonnegative and
  compartmental dynamical systems}, Princeton University Press, 2010.

\bibitem[HD96]{Hees}
JAP Heesterbeek and Klaus Dietz, \emph{The concept of ro in epidemic theory},
  Statistica neerlandica \textbf{50} (1996), no.~1, 89--110.

\bibitem[HK19]{Hurtado19}
Paul~J Hurtado and Adam~S Kirosingh, \emph{Generalizations of the ‘linear
  chain trick': incorporating more flexible dwell time distributions into mean
  field ode models}, Journal of mathematical biology \textbf{79} (2019), no.~5,
  1831--1883.

\bibitem[HT81]{hun}
Vera H{\'a}rs and J{\'a}nos T{\'o}th, \emph{On the inverse problem of reaction
  kinetics}, Qualitative theory of differential equations \textbf{30} (1981),
  363--379.

\bibitem[JKVdD77]{jeffries1977matrix}
Clark Jeffries, Victor Klee, and Pauline Van~den Driessche, \emph{When is a
  matrix sign stable?}, Canadian Journal of Mathematics \textbf{29} (1977),
  no.~2, 315--326.

\bibitem[Kal12]{kaltenbach2012unified}
Hans-Michael Kaltenbach, \emph{A unified view on bipartite species-reaction and
  interaction graphs for chemical reaction networks}, arXiv preprint
  arXiv:1210.0320 (2012).

\bibitem[LDA13]{letellier2013can}
Christophe Letellier, Fabrice Denis, and Luis~A Aguirre, \emph{What can be
  learned from a chaotic cancer model?}, Journal of theoretical biology
  \textbf{322} (2013), 7--16.

\bibitem[Lic21]{Licht}
Daniel Lichtblau, \emph{Symbolic analysis of multiple steady states in a mapk
  chemical reaction network}, Journal of Symbolic Computation \textbf{105}
  (2021), 118--144.

\bibitem[MA94]{Mir-An}
R~Miron and M~Anastasiei, \emph{The geometry of {L}agrange spaces: theory and
  applications}, Springer Science \& Business Media, 1994.

\bibitem[MC08]{mincheva2008multigraph}
Maya Mincheva and Gheorghe Craciun, \emph{Multigraph conditions for
  multistability, oscillations and pattern formation in biochemical reaction
  networks}, Proceedings of the IEEE \textbf{96} (2008), no.~8, 1281--1291.

\bibitem[Min12]{mincheva2012oscillations}
Maya Mincheva, \emph{Oscillations in non-mass action kinetics models of
  biochemical reaction networks arising from pairs of subnetworks}, Journal of
  Mathematical Chemistry \textbf{50} (2012), no.~5, 1111--1125.

\bibitem[Mir01]{Miron-Hr-Shim-Sab}
Radu Miron, \emph{The geometry of {H}amilton and {L}agrange spaces}, vol. 118,
  Springer Science \& Business Media, 2001.

\bibitem[MR07]{mincheva2007graph}
Maya Mincheva and Marc~R Roussel, \emph{Graph-theoretic methods for the
  analysis of chemical and biochemical networks. i. multistability and
  oscillations in ordinary differential equation models}, Journal of
  mathematical biology \textbf{55} (2007), no.~1, 61--86.

\bibitem[MY20]{macauley2020case}
Matthew Macauley and Nora Youngs, \emph{The case for algebraic biology: from
  research to education}, Bulletin of Mathematical Biology \textbf{82} (2020),
  1--16.

\bibitem[Nil22a]{NillI}
\emph{Symmetries and normalization in 3-compartment epidemic models. i: The
  replacement number dynamics.}, arXiv:2301.00159 (2022).

\bibitem[Nil22b]{NillII}
\emph{Symmetries and normalization in 3-compartment epidemic models. ii:
  Equilibria and stability.}, arXiv:2301.00159 (2022).

\bibitem[NL23]{Nea-Lit}
M.~Neagu and A.V. Litr\u{a}, \emph{{SIR} dynamical model with demography and
  {L}agrange-{H}amilton geometries}, Scientific Studies \& Research. Series
  Mathematics \& Informatics \textbf{33} (2023), no.~1, 87--96.

\bibitem[NO22]{Nea-Oana}
Mircea Neagu and Alexandru Oan{\u{a}}, \emph{Dual jet geometrization for
  time-dependent {H}amiltonians and applications}, Springer Nature, 2022.

\bibitem[NU01]{Nea-Udr}
M.~Neagu and C.~Udri\c{s}te, \emph{From {PDE}s systems and metrics to
  multi-time field theories and geometric dynamics}, Semin. Mech.-Differ. Dyn.
  Syst., West Univ. Timi\c{s}oara, Romania \textbf{79} (2001), 1--33.

\bibitem[OMBA12]{otero2012characterizing}
Irene Otero-Muras, Julio~R Banga, and Antonio~A Alonso, \emph{Characterizing
  multistationarity regimes in biochemical reaction networks}, PLoS One
  \textbf{7} (2012), no.~7, e39194.

\bibitem[OMYS17]{otero2017chemical}
Irene Otero-Muras, Pencho Yordanov, and Joerg Stelling, \emph{Chemical reaction
  network theory elucidates sources of multistability in interferon signaling},
  PLoS computational biology \textbf{13} (2017), no.~4, e1005454.

\bibitem[OSS22]{Ott}
Stefania Ottaviano, Mattia Sensi, and Sara Sottile, \emph{Global stability of
  sairs epidemic models}, Nonlinear Analysis: Real World Applications
  \textbf{65} (2022), 103501.

\bibitem[Ple77]{plemmons1977m}
Robert~J Plemmons, \emph{M-matrix characterizations. i—nonsingular
  m-matrices}, Linear Algebra and its Applications \textbf{18} (1977), no.~2,
  175--188.

\bibitem[PS05]{pachter2005algebraic}
Lior Pachter and Bernd Sturmfels, \emph{Algebraic statistics for computational
  biology}, vol.~13, Cambridge university press, 2005.

\bibitem[RS13]{RobSti}
Marguerite Robinson and Nikolaos~I Stilianakis, \emph{A model for the emergence
  of drug resistance in the presence of asymptomatic infections}, Mathematical
  biosciences \textbf{243} (2013), no.~2, 163--177.

\bibitem[SE22]{sadeghimanesh2022resultant}
AmirHosein Sadeghimanesh and Matthew England, \emph{Resultant tools for
  parametric polynomial systems with application to population models}, arXiv
  preprint arXiv:2201.13189 (2022).

\bibitem[TF21]{torres2021symbolic}
Ang{\'e}lica Torres and Elisenda Feliu, \emph{Symbolic proof of bistability in
  reaction networks}, SIAM Journal on Applied Dynamical Systems \textbf{20}
  (2021), no.~1, 1--37.

\bibitem[TNP18]{Toth}
J{\'a}nos T{\'o}th, Attila~L{\'a}szl{\'o} Nagy, and D{\'a}vid Papp,
  \emph{Reaction kinetics: exercises, programs and theorems}, Springer, 2018.

\bibitem[Udr00]{Udr}
C.~Udri\c{s}te, \emph{Geometric dynamics}, Kluwer Academic Publishers,
  Dordrecht, 2000.

\bibitem[VAA24]{VAA}
Nicola Vassena, Florin Avram, and Rim Adenane, \emph{Finding bifurcations in
  mathematical epidemiology via reaction network methods}, arXiv preprint
  arXiv:2405.14576 (2024).

\bibitem[VdDW02]{Van}
Pauline Van~den Driessche and James Watmough, \emph{Reproduction numbers and
  sub-threshold endemic equilibria for compartmental models of disease
  transmission}, Mathematical biosciences \textbf{180} (2002), no.~1-2, 29--48.

\bibitem[VdDW08]{Van08}
P~Van~den Driessche and James Watmough, \emph{Further notes on the basic
  reproduction number}, Mathematical epidemiology, Springer, 2008,
  pp.~159--178.

\bibitem[Wil09]{Wilh}
Thomas Wilhelm, \emph{The smallest chemical reaction system with bistability},
  BMC systems biology \textbf{3} (2009), no.~1, 1--9.

\bibitem[WRK05]{Wearing2005}
H.J Wearing, P~Rohani, and M.J. Keeling, \emph{Appropriate models for the
  management of infectious diseases}, PLoS Medicine \textbf{7} (2005), no.~2,
  621--627.

\bibitem[WSFC17]{Wang2017}
Xiaojing Wang, Yangyang Shi, Zhilan Feng, and Jingan Cui, \emph{Evaluations of
  interventions using mathematical models with exponential and non-exponential
  distributions for disease stages: the case of ebola}, Bulletin of
  mathematical biology \textbf{79} (2017), 2149--2173.

\bibitem[YB08]{YangBrauer}
Christine~K Yang and Fred Brauer, \emph{Calculation of $ r\_0 $ for
  age-of-infection models}, Mathematical Biosciences \& Engineering \textbf{5}
  (2008), no.~3, 585.

\bibitem[YCMP22]{yu2022graph}
Polly~Y Yu, Gheorghe Craciun, Maya Mincheva, and Casian Pantea, \emph{A
  graph-theoretic condition for delay stability of reaction systems}, SIAM
  Journal on Applied Dynamical Systems \textbf{21} (2022), no.~2, 1092--1118.

\end{thebibliography}

\end{document}